\documentclass[10pt]{article}
\usepackage{amsmath}
\usepackage{color}
\usepackage{fancyhdr}
\usepackage{verbatim}
\usepackage{amsfonts,longtable}
\usepackage{epsfig}
\usepackage{mathrsfs}
\usepackage{amsfonts}
\usepackage{color}

\def\R{\hbox{{\rm I}\kern-0.2em{\rm R}\kern0.2em}}
\def\a{\alpha}   \def\k{\kappa}

\def\d{{\rm d}}

\def\bn{\begin{equation}}
\def\en{\end{equation}}
\def\bny{\begin{eqnarray}}
\def\eny{\end{eqnarray}}
\def\be{\begin{eqnarray*}}
\def\ee{\end{eqnarray*}}
\def\bc{\begin{center}}
\def\ec{\end{center}}
\def\X{{\cal X}} \def\U{{\cal U}}
\def\p{\partial}
\def\({\left(}
\def\){\right  )}
\def\[{\left[}
\def\]{\right]}
\def\bc{\begin{center}}
\def\ec{\end{center}}

\def\n{\noindent}

\newtheorem{dfn}{Definition}[section]
\newtheorem{thm}{Theorem}[section]
\newtheorem{rem}{Remark}[section]
\newtheorem{pro}{Proposition}[section]

\newtheorem{cor}{Corollary}[section]
\newtheorem{lem}{Lemma}[section]
\newtheorem{exm}{Example}[section]

\def\bn{\begin{equation}}
\def\en{\end{equation}}
\def\bny{\begin{eqnarray}}
\def\eny{\end{eqnarray}}
\def\be{\begin{eqnarray*}}
\def\ee{\end{eqnarray*}}
\def\bdn{\begin{dfn}}
\def\edn{\end{dfn}}
\def\btm{\begin{thm}}
\def\etm{\end{thm}}
\def\bpf{\begin{proof}}
\def\epf{\end{proof}}
\def\bpn{\begin{pro}}
\def\epn{\end{pro}}
\def\brk{\begin{rem}}
\def\erk{\end{rem}}
\def\bcy{\begin{cor}}
\def\ecy{\end{cor}}
\def\blm{\begin{lem}}\def\elm{\end{lem}}
\def\bex{\begin{exm}}
\def\eex{\end{exm}}

\def\X{{\cal X}}  
  
 \def\R{\mathscr{R}}
\def\Q{{\cal Q}}\def\r{{\cal X}}

\begin{document}
\bc \Large{{\bf A basis of hierarchy of generalized symmetries  and their conservation laws for the (3+1)-dimensional diffusion equation }}\ec
\medskip

\bc

J J H  Bashingwa* and  A H Kara**\\
School of Mathematics, University of the Witwatersrand, \\
Johannesburg, Wits 2001, South Africa.\\\begin{footnotesize}
{* jeanjuste@aims.ac.za\\
** Abdul.Kara@wits.ac.za}
\end{footnotesize}\ec

%
%
%

\begin{abstract}
We determine, by hierarchy,  dependencies between higher order linear symmetries  which occur when generating them using recursion operators. Thus, we deduce a formula which gives  the number of independent generalized symmetries (basis) of several orders. We construct a basis for conservation laws (with respect to the group  admitted by the system of differential equation) and hence  generate infinite conservation laws in each equivalence class.
\end{abstract}

\section{Introduction}
Generalized symmetries have been shown to be important in the study of nonlinear wave equations. The possession of infinite number of such symmetries 
has a bearing on integrability. Completely integrable nonlinear partial differential equations (PDEs) have a rich structure. For example, they have infinitely many conservation laws of increasing order. These PDEs have a Lax pair, they can also be solved with the inverse scattering transform and admit soliton solutions of any order \cite{MJ1,MJ2}.    Supersymmetries are interpretable as generalized symmetries \cite{Cast}.

There are two approaches for finding generalized symmetries. The first involves the usage of the invariance condition. This method is essentially the same as the one used to find Lie symmetries although the intervening calculations usually are far more complicated as we are adding  derivatives of the relevant dependent variables in the infinitesimals. A second approach is the use of recursion operators which will generate infinite families of symmetries at once.

Let \bn \Delta (x,u_{(\a)})=0 \label{systm1}\en be a system of differential equations with $m$ dependent variables $u=(u^1,...,u^m)$  and $n$ independent variables $x=(x_1,...,x_n)$. The notation $u_{(\a)}$ stands for all derivatives of $u$ with respect to $x$ with order no greater than $\a$.

The following  can be found in {\cite{Olver93,{ka},{ka1},blu}} 
\bdn  
Consider a system of differential equations given by \eqref{systm1}.
A recursion operator for $\Delta$ is a linear operator $\mathscr{R}$ in the space of differential functions with the property that whenever $V_Q$ is an evolutionary symmetry of $\Delta$, so is $V_{\tilde{Q}}$ with $\tilde{Q}=\mathscr{R} Q$\edn

 Alot of work appears in the literature on recursion operators, symmetries and differential equations \cite{MJ3,MJ4,MJ5}. In \cite{MJ6}, the author provided an approach for constructing all recursion operators for a given linear system of differential equation.
 
 \bpn  
 Consider the system \eqref{systm1}, with $\Delta$ denoting a linear differential operator. A second linear operator $\mathscr{R}$ not depending on u or its derivatives is a recursion operator for $\Delta$ if and only if $Q=\mathscr{R}[u]$ is a characteristic of a linear generalized symmetry to the system. \epn 

Even if it becomes  easy to compute all generalized symmetries once the recursion operators are known for a system of differential equations, there are dependencies among the resulting symmetries stemming from the relation among operators. For example, the two-dimensional wave equation $$u_{tt}=u_{xx}+u_{yy}$$ has ten recursion operators \cite{Olver93}. In \cite{Del1,Del2}, the authors proved that there are  $$\frac{(2k+1)(2k+2)(2k+2)}{6}$$   $k$-th order independent symmetries generated by these recursion operators.
\bdn A conservation law of the system \eqref{systm1} is a divergence expression \text{Div} $T=D_iT^i$  which vanishes for all solutions of \eqref{systm1} i.e \edn
Div  $T|_{\Delta=0}=0$.
The n-tuple $T=(T^1,...,T^n)$  is called a conserved vector of this conservation law,  $D_i=D_{x_i}$ is the total differentiation with respect to the variable $x_i$, and $T^i$ are differential functions.

We  have the  conservation laws in differential forms as follow 
\bdn  The differential form $(n-1)$-form\edn
$\omega = T^i(x,u_{(r)})\p_{x_i} \lrcorner (\d x_1 \wedge ...\wedge \d x_n)  $ is called conserved form of \eqref{systm1} if $$D\omega  =0$$ is satisfied on the manifold in the space of variables defined by the system \eqref{systm1}

\bdn Two conserved vectors $T$ and $T'$ are called equivalent if the vector function $T-T'$ is a trivial conserved vector\label{defnn1}\edn

\bpn Any point transformation $g$ maps a class of equations in the conserved form into itself. \epn Let $g:\bar{x}=x_g(x,u), \quad \bar{u}=u_g(x,u)$ be a point transformation of the system \eqref{systm1}. The prolongation to the jet space of $g$ transforms the vector in conserved form as follow

$T^i_g(\bar{x},\bar{u}_{(r)})=\frac{D_{x_j}\bar{x_i}}{|D_x\bar{x}|}T^j(x,u_{(r)}),$ or 
\bn \bar{T}\equiv  T_g(\bar{x},\bar{u}_{(r)})=\frac{1}{|D_x\bar{x}|}(D_x\bar{x})T(x,u_{(r)})\label{conser}\en

\bdn Two conservation laws $T$ and $\bar{T}$ are equivalent with respect to the group  of transformation $G$ if there exist a transformation $g\in G$ such that $T$ and $\bar{T}\equiv T_g$ are equivalent in the sense of Definition 1.4 \edn

Once a conservation law is known, we can generate new conservation laws using the infinitesimal generators $X=\xi ^i\p_i + \eta ^ \mu \p_{u^\mu}$ corresponding to the system \eqref{systm1}.  In fact differentiating \eqref{conser} with respect to the parameter $\epsilon$ and taking the value $\epsilon=0$ we get see ( \cite{pop})

 \bn \bar{T^i}= -X_{}T^i+ (D_j\xi ^i)T^j -(D_j \xi ^j)T^i \label{tbar}\en  

A similar formula for symmetry in characteristic form ($\xi^i=0$) was known  earlier \cite{Olver93, ibra}. It has been used by various authors \cite{ka, kam} to find basis of conservation law, i.e a set which generates a whole set of conservation law by action of generalized symmetry operator and operation of linear combination.

\btm  Let  ad $Y(X)=Z$ (ad $Y(X)=[X,Y].$)such that $Y$ is associated with the conserved vector $T$ of \eqref{systm1}. Then $\bar{T}$ given by \eqref{tbar} is trivial conserved vector if $Z=b Y$, for any constant b.\etm 

\btm  Suppose that $$X=\xi ^i\p_{x^i}+ \eta ^\a \p_{u^\a}+\zeta^\a_i\p_{u^\a_i}+\zeta^\a_{i_1i_2}\p_{u^\a_{
i_1i_2}}+...$$  where the coefficient are given by

$\zeta^\a_i=D_i (\eta^\a)-u_j^\a D_i(\xi^j)\\
\zeta^\a_{i_1...i_s}=D_{i_s}(\zeta^\a_{i_1...i_{s-1}})-u^\a_{ji_1...i_{s-1}}D_{i_s}(\xi^j),\quad s>1,$

is a Lie-B\"acklund operator of the system \eqref{systm1}
such that the  conserved form $\omega$ of a system  \eqref{systm1} is invariant under $X$. Then 
\bn  X_{}T^i- (D_j\xi ^i)T^j +(D_j \xi ^j)T^i=0 \label{tbr}\en  

 Here $D_i$ is the total derivative operator.\etm

%
%

In this work, we are interested in the
 3-dimensional diffusion equation

\bn U_t=U_{xx}+U_{yy}+U_{zz}. \label{heat}\en

In the second section, we compute the Lie-point symmetry generators of \eqref{heat}. By  Proposition 1.1, we can deduce the recursion operators for \eqref{heat}. We construct dependencies which occur among operators by hierarchy meaning that we start by first order, second, third, ... and dependencies in $k$-th order are constructed using dependencies in $(k-1)$-th order. We generalize our method to any order $n$ and find the number of independent symmetries generated by these recursion operators.

In the third section, even though it is cumbersome, we construct the basis of all generalized symmetries (first, second and third order) using the invariance condition and compare the results with those in the first section.

In the fourth section, we construct the basis of conservation laws for the equation \eqref{heat}. By applying the generalized symmetries found in previous Section, we can generate infinite conservation laws.

\section{Generalized symmetries using recursion operators.}

In this Section, we construct the basis for the generalized symmetries of several orders using recursion operators.
\subsection{First order symmetries}

The Lie algebra of point symmetry generators with basis (Lie symmetries) of \eqref{heat} is given by

\bn
\begin{array}{llll}
\X_1&= 2t\p_x-Ux\p_U\\
\X_2&=2t\p_y-Uy\p_U\\
\X_3&=2t\p_z-Uz\p_U\\
\X_4& = -x\p_y+y\p_x \\
\X_5& = -x\p_z +z\p_y \\
\X_6& = \p_x \\
\X_7& =  -y\p_z+z\p_y\\
\X_8& =  \p_y\\
\X_9& =  \p_z\\
\X_{10}& =  U\p_U\\
\X_{11}& =\p_t  \\
\X_{12}& = 2t \p_t +x\p_x + y\p_y +z\p_z  \\
\X_{13}& =4t^2 \p_t +4tx\p_x+4ty\p_y+4tz\p_z-u(6t +x^2 +y^2+z^2)\p_U  \\
\label{gen0}  \\
\end{array}
\en
and the indefinite generator of symmetry
$$\X_{14} =f_{1}(x,y,z,t)\p_U$$ where the function $f_1$ satisfies $f_{1,t}-f_{1,xx}-f_{1,yy}f_{1,zz}=0$

\eqref{gen0} forms a Lie algebra with the following commutators
\bn
\begin{array}{lllll}
&[\r_1, \r_4]=-\r_2, & [\r_2,\r_8]=\r_{10},& [\r_4,\r_5]=\r_7, & [\r_5,\r_6]=\r_9, \\

&[\r_1, \r_5]=-\r_3, & [\r_3,\r_5]=\r_{1},& [\r_4,\r_6]=\r_8,& [\r_5,\r_7]=\r_4, \\

&[\r_1, \r_6]=\r_{10},&[\r_3,\r_7]=\r_{2},&[\r_4,\r_7]=-\r_5,& [\r_5,\r_9]=-\r_6 ,\\
&[\r_1,\r_{11}]=-2\r_6, & [\r_2,\r_{11}]=-2\r_8,&[\r_3,\r_{11}]=-2\r_9,&[\r_6,\r_{12}]=\r_6,\\

&[\r_1,\r_{12}]=-\r_1, & [\r_2,\r_{12}]=-\r_2,&[\r_3,\r_{12}]=-\r_3, &[\r_6,\r_{13}]=2\r_1,\\

&[\r_2, \r_4]=\r_{1},& [\r_3,\r_9]=\r_{10},&[\r_4,\r_8]=-\r_6,& [\r_7,\r_8]=\r_9 ,\\

&[\r_2, \r_7]=-\r_{3},& [\r_8,\r_{12}]=\r_8,&[\r_8,\r_{13}]=2\r_2,& [\r_7,\r_9]=-\r_8 ,\\

&[\r_9, \r_{12}]=\r_{9},& [\r_9,\r_{13}]=2\r_3,&[\r_{11},\r_{12}]=2\r_{11},& [\r_{11},\r_{13}]=4\r_{12}-6\r_{10} ,\\
&[\r_{12},\r_{13}]=2\r_{13}.
\label{alge}\end{array}
\en 

all the remaining commutators are zeros. It is important to notice that the first ten generators $\X_1,...,\X_{10}$ form a sub-algebra of the Lie algebra \eqref{alge} which is the following

\bn
\begin{array}{lllll}
&[\r_1, \r_4]=-\r_2, & [\r_2,\r_8]=\r_{10},& [\r_4,\r_5]=\r_7, & [\r_5,\r_6]=\r_9, \\

&[\r_1, \r_5]=-\r_3, & [\r_3,\r_5]=\r_{1},& [\r_4,\r_6]=\r_8,& [\r_5,\r_7]=\r_4, \\

&[\r_1, \r_6]=\r_{10},&[\r_3,\r_7]=\r_{2},&[\r_4,\r_7]=-\r_5,& [\r_5,\r_9]=-\r_6 ,\\

&[\r_2, \r_4]=\r_{1},& [\r_3,\r_9]=\r_{10},&[\r_4,\r_8]=-\r_6,& [\r_7,\r_8]=\r_9 ,\\

&[\r_2, \r_7]=-\r_{3},& && [\r_7,\r_9]=-\r_8 ,\\

\label{alge1}\end{array}
\en 

Among all the generalized symmetries vector fields, we are interested in evolutionary vector field which are given by

$$V_\Q=\Q [U]\p_U$$

where [U] is includes the basis, fibre and jet variables that the characteristic $Q$ depends on.

From \eqref{gen0}, and with the usage of Proposition 1.1, we see that they are nine recursion operators, namely

\bn
\begin{array}{ll}
\R_1=& 2t \p_x +x\\\R_2= &2t\p_y+y\\\R_3=&2t\p_z+z\\\R_4=&-x\p_y+y\p_x\\\R_5=&-x\p_z+z\p_x\\\R_6=&\p_x\\\R_7=&-y\p_z+z\p_y\\\R_8=&\p_y\\\R_9=&\p_z
\label{rec}\end{array}
\en

Hence,  every first order generalized symmetry has as its characteristic $\Q$ a
linear constant coefficient combination of the following ten independent "basic" characteristic 

\bn \begin{array}{ll}
\Q_{0}=& U\\
\Q_1^1= & 2t U_x +xU\\ \Q_2^1=& 2t U_y +yU \\
\Q_3^1=& 2t U_z + zU\\ \Q_4^1= & -xU_y + yU_x \\ \Q_5^1=& -x U_z + zU_y\\ \Q_6^1= &U_x \\ \Q_7^1= & -y U_z + zU_y\\ \Q_8^1 =& U_y\\ \Q_9^1=&U_z\\ \label{gen1}
\end{array}
\en

obtained by applying successively the recursion operators in \eqref{rec} on $\Q_0=U$; i.e $$\Q_1^1=\R_1[U],\quad \Q_2^1=\R_2[U],...,\Q_9^1=\R_9[U]$$ 

Note*  : If there is no confusion we will omit $[U]$ in our notation and use $\R_i$ instead of $\R_i[U]$ for these characteristics.

\subsection{Second order symmetries}
We need to construct the basis "characteristics" of all second order generalized symmetries using the recursion operators. Basically it is comprised of  the ten first order "basic characteristic" plus all second order combination with repetitions $\R_i\R_j$  with $i,j=1...9$ and $i\leq j$

Note that we don't need to compute the cases $i>j$ because of the relations in \eqref{alge1}

The total number of the all second order combinations with repetition obtained with 9 recursion operators is 45 i.e $$\left( \begin{matrix}
9+2-1 \\ 
2
\end{matrix}\right)=\frac{10!}{8!2!} $$ 

It is turns out that some of the first order characteristics can be expressed in terms of the second order ones as follow  

\bn
\begin{array}{ll}
\R_4=& -\R_1\R_8+\R_2\R_6\\
\R_5=&-\R_1\R_9+\R_3\R_6\\
\R_7=& -\R_2\R_9+\R_3\R_8
\end{array}\label{fir-od}
\en
There is also two dependencies between symmetries generated by the second order combinations of the recursion operators: 
\bn
\begin{array}{ll}
\R_3\R_4=&\R_2\R_5-\R_1\R_7\\
\R_6\R_7=&\R_5\R_8-\R_4\R_9

\end{array}\label{sec-od}
\en

We need to take out these dependencies. Hence, every second order generalized symmetry has as its characteristic $\Q$ a
linear constant coefficient combination of ($45+10-5=50$) fifty "basic" characteristic 
\bn
\begin{array}{llllllllll}
\R_1,&\R_1\R_1, \\
\R_2,&\R_1\R_2,&\R_2\R_2,\\
\R_3,&\R_1\R_3,&\R_2\R_3,&\R_3\R_3,\\
&\R_1\R_4,&\R_2\R_4,&&\R_4\R_4,\\
&\R_1\R_5,&\R_2\R_5,&\R_3\R_5,&\R_4\R_5,&\R_5\R_5,\\
\R_6,&\R_1\R_6,&\R_2\R_6,&\R_3\R_6,&\R_4\R_6,&\R_5\R_6,
&\R_6\R_6
,\\
&\R_1\R_7,&\R_2\R_7,&\R_3\R_7,&\R_4\R_7,&\R_5\R_7,&
&\R_7\R_7,\\
\R_8,&\R_1\R_8,&\R_2\R_8,&\R_3\R_8,&\R_4\R_8,&\R_5\R_8,&
\R_6\R_8
,&\R_7\R_8,&\R_8\R_8,\\
\R_9,&\R_1\R_9,&\R_2\R_9,&\R_3\R_9,&\R_4\R_9,&\R_5\R_9,&
\R_6\R_9
,&\R_7\R_9,&\R_8\R_9,&\R_9\R_9.\\
\end{array}
\en

plus the characteristic $\Q_0$

\subsection{Third order symmetries}

The total number of all third order combinations with repetition  $R_iR_jR_k$ where ($i\leq j\leq k$ and $i,j,k=1,...,9$) is 165 i.e
$$\left(\begin{matrix}
9+3-1\\3
\end{matrix} \right)=\frac{11!}{8!3!}$$ 
Basically, the basis of all third order generalized symmetries should include  the "fifty" second order symmetries plus all third order combinations with repetition .

Now let  construct dependencies by hierarchy. 
We see from \eqref{fir-od}, some of second order characteristics can be expressed in terms third order ones or combinations of second and third order. We got the following 22 dependencies

\bn
\begin{array}{llllll}
\R_1\R_4=&-\R_1\R_1\R_8+\R_1\R_2\R_6,\quad &
\R_4\R_4=&-\R_1\R_4\R_8+\R_2\R_4\R_6-\R_2\R_8-\R_1\R_6\\
\R_2\R_4=&-\R_1\R_2\R_8+\R_2\R_2\R_6,\quad &

\R_4\R_5=&-\R_1\R_5\R_8+\R_2\R_5\R_6-\R_2\R_9\\
\R_4\R_6=&-\R_1\R_6\R_8+\R_2\R_6\R_6,\quad&
\R_4\R_8=&-\R_1\R_8\R_8+\R_2\R_6\R_8\\
\R_4\R_9=&-\R_1\R_8\R_9+\R_2\R_6\R_9,\quad &
\R_1\R_5=& -\R_1\R_1\R_9+\R_1\R_3\R_6 \\
\R_2 \R_5 = & -\R_1 \R_2 \R_9 +\R_2 \R_3 \R_6,\quad &
\R_3\R_5=& -\R_1\R_3\R_9+\R_3\R_3\R_6\\
\R_5\R_5=& -\R_1\R_5\R_9-\R_3\R_9+\R_3\R_5\R_6-\R_1\R_6,\quad &\R_5\R_6=&-\R_1\R_6\R_9+\R_3\R_6\R_6\\
\R_5\R_7=&-\R_1\R_7\R_9-\R_1\R_8+\R_3\R_6\R_7,&
\R_5\R_8=&-\R_1\R_8\R_9+\R_3\R_6\R_6\\
\R_5\R_9=&-\R_1\R_9\R_9+\R_3\R_6\R_9,&
\R_1\R_7=&-\R_1\R_2\R_9+\R_1\R_3\R_8\\
\R_2\R_7=&-\R_2\R_2\R_9+\R_2\R_3\R_8,&
\R_3\R_7=&-\R_2\R_3\R_9+\R_3\R_3\R_8\\
\R_4\R_7=&-\R_2\R_4\R_9+\R_1\R_9+\R_3\R_4\R_8,&
\R_7\R_7=&-\R_2\R_7\R_9+\R_3\R_7\R_8-\R_3\R_9-\R_2\R_8\\
\R_7\R_8=&-\R_2\R_8\R_9+\R_3\R_8\R_8,&
\R_7\R_9=&-\R_2\R_9\R_9+\R_3\R_8\R_9\\

\end{array}\label{thi_1}
\en

and from \eqref{sec-od} we got the following 18 dependencies between the third order characteristics

\bn
\begin{array}{lllll}
\R_1\R_3\R_4=&\R_1\R_2\R_5-\R_1\R_1\R_7,&\R_3\R_4\R_4
=&\R_2\R_4\R_5+\R_2\R_2\R_9-\R_2\R_3\R_8-

\R_1\R_4\R_7\\&&& +\R_1\R_1\R_9+\R_1\R_3\R_6\\\R_2\R_3\R_4
=&\R_2\R_2\R_5-\R_1\R_2\R_7,&\R_3\R_4\R_6=&\R_2\R_5\R_6-\R_1\R_5\R_8+\R_1\R_4\R_9\\\R_3\R_3\R_4=&\R_2\R_3\R_5-\R_1\R_3\R_7,&

\R_3\R_4\R_5=&\R_2\R_5\R_5-\R_1\R_5\R_7-\R_1\R_1\R_8+\R_1\R_2\R_6\\

\R_3\R_4\R_7=&\R_2\R_5\R_7-\R_1\R_7\R_7,&
\R_3\R_4\R_8=&\R_2\R_5\R_8-\R_1\R_7\R_8\\
\R_3\R_4\R_9=&\R_2\R_5\R_9-\R_1\R_7\R_9,&
\R_1\R_6\R_7=& \R_1\R_5\R_8-\R_1\R_4\R_9\\
\R_2\R_6\R_7=& \R_2\R_5\R_8-\R_2\R_4\R_9,\quad&
\R_3\R_6\R_7=& \R_3\R_5\R_8-\R_2\R_5\R_9+\R_1\R_7\R_9\\
\R_4\R_6\R_7=& \R_4\R_5\R_8-\R_4\R_4\R_9,&
\R_5\R_6\R_7=& \R_5\R_5\R_8-\R_4\R_5\R_9-\R_2\R_9\R_9+\R_3\R_8\R_9\\
\R_6\R_6\R_7=& \R_5\R_6\R_8-\R_4\R_6\R_9,&
\R_6\R_7\R_7=& \R_5\R_7\R_8-\R_1\R_9\R_9-\R_3\R_6\R_9-\R_4\R_7\R_9\\&&&-\R_1\R_8\R_8-\R_2\R_6\R_9\\
\R_6\R_7\R_8=& \R_5\R_8\R_8-\R_4\R_8\R_9,&
\R_6\R_7\R_9=& \R_5\R_8\R_9-\R_4\R_9\R_9\\
\end{array}\label{thi_2}
\en

Total number of dependencies is 18+22=40. We will deduce this number to get the exact number of non independent third order characteristic.

Hence, every third order generalized symmetry has its characteristic which is  a
linear constant coefficient combination of ($165+50-40=175$) hundred and seventy-five "basic" characteristic.

\subsection{Fourth order symmetries}

The total number of all fourth order combinations with repetition  constructed with nine recursion operators is 495 i.e

$$\left(\begin{matrix}
9+4-1\\4
\end{matrix} \right)$$ 

The basis of the fourth order generalized symmetries is included those fourth order combinations plus the 175 third order symmetries.

Let now construct dependencies by hierarchy.
From relations \eqref{thi_1} we got 91 dependencies , and from \eqref{thi_2} we got 89. We will deduce those dependencies to obtain the exact number for the basis.

Hence, every fourth order generalized symmetry has  characteristic which is  a
linear constant coefficient combination of ($495+175-(91+89)=490$) four hundred and ninety. "basic" characteristic.

\subsection{Generalization}

We have seen that, for a given order $n$, there are two kinds of dependencies:

- between $n$-th order symmetries ,

-  between $n$-th order and $(n-1)$-th order symmetries.

$\bullet$For the second order, we have 5 dependencies

2 and 3

$\bullet$ For the third order, we have  40 dependencies
$$18=2\left(\begin{matrix}
9+1-1\\1
\end{matrix} \right)$$ and $$22=3\left(\begin{matrix}
9+1-1\\1
\end{matrix} \right)-(2(1+0)+2(1-0)+1)$$

$\bullet$ For the fourth order, we have 180 dependencies

$$89=2\left(\begin{matrix}
9+2-1\\2
\end{matrix} \right)-1$$ and 

$$91=3\left(\begin{matrix}
9+2-1\\2
\end{matrix} \right)-(2(9+1)+1+1+2(9-(1+1))+9-1)$$

$\bullet$ For the fifth order, we have 601 dependencies

$$321=2\left(\begin{matrix}
9+3-1\\3
\end{matrix} \right)-9$$ and $$280=3\left(\begin{matrix}
9+3-1\\3
\end{matrix} \right)-2(45+(9-1)+9+9+2(45-(9+8))+(45-(9+1)))$$

$\bullet$ For the sixth order, we have 1659 dependencies

$$945=22\left(\begin{matrix}
9+4-1\\4
\end{matrix} \right)-45$$ and 

$$714=32\left(\begin{matrix}
9+4-1\\4
\end{matrix} \right)-(2(165+(45-9))+45+(45-1)+2(165-(45+(45-9)))+165
-(45+(9-1)))$$
For each order $k$ we have

\bn
\begin{array}{lll}
2\left(\begin{matrix}
9+k-2-1\\k-2
\end{matrix} \right)-\left(\begin{matrix}
9+k-4-1\\k-4
\end{matrix} \right)&=2\left(\begin{matrix}
k+6\\k-2
\end{matrix} \right)-\left(\begin{matrix}
k+4\\k-4
\end{matrix} \right)\\&=\left(\begin{matrix}
k+7\\k-1
\end{matrix} \right)+\left(\begin{matrix}
k+6\\k-1
\end{matrix} \right)-3\left(\begin{matrix}
k+5\\k-1
\end{matrix} \right)+\left(\begin{matrix}
k+4\\k-1
\end{matrix} \right)\\&=
\left(\begin{matrix}
k+8\\k
\end{matrix} \right)-4\left(\begin{matrix}
k+6\\k
\end{matrix} \right)+4\left(\begin{matrix}
k+5\\k
\end{matrix} \right)-\left(\begin{matrix}
k+4\\k
\end{matrix} \right)

\end{array}
\en

dependencies between $k$-th order symmetries generated by the 9 recursion operators and 

\bn
\begin{array}{llll}
3\left(\begin{matrix}
k+7\\k-1
\end{matrix} \right)-5\left(\begin{matrix}
k+6\\k-2
\end{matrix} \right)+\left(\begin{matrix}
k+5\\k-3
\end{matrix} \right)+\left(\begin{matrix}
k+4\\k-4
\end{matrix} \right)&=4\left(\begin{matrix}
k+5\\k-1
\end{matrix} \right)-\left(\begin{matrix}
k+4\\k-1
\end{matrix} \right)\\&=

4\left(\begin{matrix}
k+6\\k
\end{matrix} \right)-5\left(\begin{matrix}
k+5\\k
\end{matrix} \right)+\left(\begin{matrix}
k+4\\k
\end{matrix} \right)
\end{array}
\en

dependencies between $(k-1)$  and $k$ order symmetries

For a given order $n$, the  number $N$ of independent symmetries generated by the recursion operator is

\bn
\begin{array}{lll}
N&=\sum\limits_{k=0}^{n} \left(\begin{matrix}
9+k-1\\k
\end{matrix} \right)-\sum\limits _{k=0}^{n-1}\[4\left(\begin{matrix}
k+6\\k
\end{matrix} \right)-5\left(\begin{matrix}
k+5\\k
\end{matrix} \right)+\left(\begin{matrix}
k+4\\k
\end{matrix} \right)\]\\&-\sum\limits^{n}_{k=0}\[\left(\begin{matrix}
k+8\\k
\end{matrix} \right)-4\left(\begin{matrix}
k+6\\k
\end{matrix} \right)+4\left(\begin{matrix}
k+5\\k
\end{matrix} \right)-\left(\begin{matrix}
k+4\\k
\end{matrix} \right)
\]\\
&=-\sum\limits _{k=0}^{n-1}\[4\left(\begin{matrix}
k+6\\k
\end{matrix} \right)-5\left(\begin{matrix}
k+5\\k
\end{matrix} \right)+\left(\begin{matrix}
k+4\\k
\end{matrix} \right)\]\\&-\sum\limits^{n}_{k=0}\[-4\left(\begin{matrix}
k+6\\k
\end{matrix} \right)+4\left(\begin{matrix}
k+5\\k
\end{matrix} \right)-\left(\begin{matrix}
k+4\\k
\end{matrix} \right)
\]\label{toata}
\end{array}
\en
We know that 

$$\sum\limits_{k=0}^{n} \left(\begin{matrix}
k+r\\k
\end{matrix} \right)=\left(\begin{matrix}
n+r+1\\r+1
\end{matrix} \right)$$

\eqref{toata} becomes

\bn
\begin{array}{lll}
N=&-4\left(\begin{matrix}
n+6\\7
\end{matrix} \right)+5\left(\begin{matrix}
n+5\\6
\end{matrix} \right)-\left(\begin{matrix}
n+4\\5
\end{matrix} \right)+4\left(\begin{matrix}
n+7\\7
\end{matrix} \right)-4\left(\begin{matrix}
n+6\\6
\end{matrix} \right)+\left(\begin{matrix}
n+5\\5
\end{matrix} \right)\\
&=5\left(\begin{matrix}
n+5\\n-1
\end{matrix} \right)+\left(\begin{matrix}
n+4\\n
\end{matrix} \right)\\&=\frac{(n+4)(n+3)^2(n+2)^2(n+1)}{144}
\end{array}
\en

They are $\frac{(n+4)(n+3)^2(n+2)^2(n+1)}{144}$ independents $n$-th order symmetries generated by these recursion operators.

For instance,

$n=1\Rightarrow N=10\\n=2\Rightarrow N=50\\n=3\Rightarrow N=175\\n=4\Rightarrow N=490$ 

In fact there 10 independent first order, 50 independent second order, 175 independent third order and 490  fourth order independent generalized symmetries.

\section{Examples, using invariance condition}

In this Section, we construct the basis of the generalized symmetries using the invariance condition.

\subsection{First order}

We let the characteristic depending on independent, dependent variables and their first order derivatives.
$\Q=\Q(x, y, z, t, U, U_x , U_y , U_z)$

we dont include $U_t$  to exclude dependency. Wherever it appears it is
replaced by $$U_{xx}+U_{yy}+U_{zz}$$ ie

$U_t=U_{xx}+U_{yy}+U_{zz}\\U_{tx}=U_{xxx}+U_{xyy}+U_{xzz}\\
...$

The  invariance condition is

\bn D_t\Q=D_x^2\Q +D_y^2\Q+D_z^2\Q\label{inv}\en

where the total derivative $D_x,D_y,D_z$ and $D_t$ are given by

\bn
\begin{array}{lll}
&D_t=\frac{\p}{\p t}+ u_t\frac{\p}{\p u}+u_{tt}\frac{\p}{\p u_t}+u_{tx}\frac{\p}{\p u_x}+u_{ty}\frac{\p}{\p u_y}+u_{tz}\frac{\p}{\p u_z}+u_{ttt}\frac{\p}{\p u_{tt}}+...,\\
&D_x=\frac{\p}{\p x}+ u_x\frac{\p}{\p u}+u_{tx}\frac{\p}{\p u_t}+u_{xx}\frac{\p}{\p u_x}+u_{xy}\frac{\p}{\p u_y}+u_{xz}\frac{\p}{\p u_z}+u_{xtt}\frac{\p}{\p u_{tt}}+...,\\
&D_y=\frac{\p}{\p y}+ u_y\frac{\p}{\p u}+u_{yt}\frac{\p}{\p u_t}+u_{yx}\frac{\p}{\p u_x}+u_{yy}\frac{\p}{\p u_y}+u_{yz}\frac{\p}{\p u_z}+u_{ytt}\frac{\p}{\p u_{tt}}+...,\\
&D_z=\frac{\p}{\p z}+ u_z\frac{\p}{\p u}+u_{zt}\frac{\p}{\p u_t}+u_{zx}\frac{\p}{\p u_x}+u_{zy}\frac{\p}{\p u_y}+u_{zz}\frac{\p}{\p u_z}+u_{ztt}\frac{\p}{\p u_{tt}}+...,\\
\end{array}
\en

Solving the invariance condition \eqref{inv}, we got the following set of first order characteristics

\bn \begin{array}{ll}
\Q_1^1= & 2t U_x +xU\\ \Q_2^1=& 2t U_y +yU \\
\Q_3^1=& 2t U_z + zU\\ \Q_4^1= & -xU_y + yU_x \\ \Q_5^1=& -x U_z + zU_y\\ \Q_6^1= &U_x \\ \Q_7^1= & -y U_z + zU_y\\ \Q_8^1 =& U_y\\ \Q_9^1=&U_z\\ \Q_{10}^1=& U\label{gen1}
\end{array}
\en
%
%
%
%
%
%

\subsection{Second order}
We let the characteristic depending up to second order derivatives, 

$Q(x,y,z,t,U,U_x,U_y,U_z,U_{xx},U_{xy},U_{xz},U_{yy},U_{yz},
,U_{zz})$

Solving the invariance condition \eqref{inv}, we got
\bn
\begin{array}{lll}
\Q_1^2 = & 4t^2 U_{xx}+ 4tx U_x +2tU +x^2U &= \R_1\R_1\\
\Q_2^2= & 4t^2(U_{yy}- U_{xx})+ 4t(yU_y-x U_x) +(y^2- x^2)U&=-\R_1\R_1+\R_2\R_2\\
\Q_3^2=& 4t^2 U_{xy}+2tyU_x+x(2t U_y +y U) &=\R_1\R_2\\
\Q_4^2=&4t^2 U_{yz}+2tz U_y +y(2tU_z +zU) &=\R_2\R_3\\
\Q_5^2= & 2tyU_{zz}-2tzU_{yz}-z(zU_y-yU_z)+yU &=\R_2-\R_3\R_7\\
\Q_6^2=& 4t^2(U_{zz}-U_{xx})+4t(zU_z-xU_x)+U(z^2-x^2)&=-\R_1\R_1+\R_3\R_3 \\
\Q_7^2= & 2txU_{yy}-2tyU_{xy}-y(yU_x-xU_y)+xU &=\R_1-\R_2\R_4\\
\Q_8^2= & U &=\R_{10}\\
\Q_9^2=& 2t^2U_{xz}+2tzU_x + x( 2tU_z + z U)&=\R_1\R_3\\
\Q_{10}^2=& 2t U_z +zU&=R_3 \\
\Q_{11}^2=& -2txU_{yy}+2tU_x + 2tyU_{xy}+y(yU_x-xU_y)&=\R_2\R_4\\
\Q_{12}^2=&2t(yU_{xx}-yU_{zz}-xU_{xy}+zU_{yz})-x(xU_y-yU_x)-z(yU_z-zU_y)&=\R_1\R_4+\R_3\R_7\\
\Q_{13}^2=&-2tU_z-2txU_{xz}+2tzU_{xx}+xzU_x-x^2U_z&=\R_1\R_5\\
\Q_{14}^2=&x(U_x-xU_{yy}+2yU_{xy})-y(yU_{xx}-yU_{zz}+2zU_{yz})+z^2U_{yy}+2tU_{zz}&=\R_7\R_7-\R_4\R_4+
\R_3\R_9\\
\Q_{15}^2=&2t(xU_{yy}-xU_{zz}-yU_{xy}+zU_{xz})+x(yU_y-zU_z)-y^2U_x+z^2U_x&=\R_3\R_5-\R_2\R_4\\
\Q_{16}^2=&-2txU_{yz}+2tzU_{xy}+y(zU_x-xU_z)&=\R_2\R_5\\
\Q_{17}^2=&2tU_{xz}+zU_x&=\R_3\R_6\\
\Q_{18}^2=&-xyU_{zz}+xzU_{yz}+yzU_{xz}-z^2U_{xy}+yU_x&=-\R_5\R_7-\R_1\R_8+\R_2\R_6\\
\Q_{19}^2=& U_x&=\R_6\\
\Q_{20}^2=&-2tyU_{zz}+2tU_y+2tz U_{yz}+z(zU_y-yU_z)&=\R_3\R_7\\
\Q_{21}^2=&-2tU_z -2tyU_{yz}+2tzU_{yy}+y(zU_y-yU_z)&=\R_2\R_7\\
\Q_{22}^2=&-y^2U_{zz}+2yzU_{yz}-z^2U_{yy}-2tU_{zz}+yU_y&=-\R_7\R_7-\R_3\R_9\\
\Q_{23}^2=&-2tyU_{xz}+2tzU_{xy}-xyU_{z}+xzU_y&=\R_1\R_7\\
\Q_{24}^2=&xyU_{zz}-xU_y-xzU_{yz}-yzU_{xz}+z^2U_{xy}&=\R_5\R_7\\
\Q_{25}^2=&2tU_{yz}+zU_y&=\R_3\R_8\\
\Q_{26}^2=&U_y&=\R_8\\
\Q_{27}^2=&2tU_{zz}+zU_z&=\R_3\R_9\\
\Q_{28}^2=&2tU_{xz}+xU_z&=\R_1\R_9\\
\Q_{29}^2=&2tU_{yz}+yU_z&=\R_2\R_9\\
\Q_{30}^2=&U_z&=\R_9\\
\end{array}
\en

\bn
\begin{array}{llll}
\Q_{31}^2&=&x(xU_{yy}-2yU_{xy})+(y^2+2t)U_{xx}-U_{zz}(y^2+2t)+2yzU_{yz}-z^2U_{yy}\\&=&\R_4\R_4-\R_7\R_7+\R_1\R_6-\R_3\R_9\\
\Q_{32}^2&=&U_{yy}(z^2-x^2)+U_{zz}(x^2+y^2+4t)+2xyU_{xy}-2xzU_{xz}+U_{xx}(z^2-y^2)-2yzU_{yz}\\&=&\R_3\R_9-\R_4\R_4+\R_7\R_7+\R_5\R_5\\
\Q_{33}^2&=&x^2U_{yz}-xyU_{xz}-xzU_{xy}+yzU_{xx}+2tU_{yz}&=\R_4\R_5+\R_2\R_9\\
\Q_{34}^2&=&-xU_{xz}+zU_{xx}&=\R_5\R_6\\
\Q_{35}^2&=&-xU_{xy}+yU_{xx}&=\R_4\R_6\\
\Q_{36}^2&=&U_{xx}&=\R_6\R_6\\
\Q_{37}^2&=&xyU_{yz}-xzU_{yy}-y^2U_{xz}+yzU_{xy}-2tU_{xz}&=\R_4\R_7-\R_1\R_9\\
\Q_{38}^2&=&-xU_{yy}+yU_{xy}&=\R_4\R_8\\
\Q_{39}^2&=&xyU_{zz}-xzU_{yz}-yzU_{xz}+z^2U_{xy}+2tU_{xy}&=\R_5\R_7+\R_1\R_8\\
\Q_{40}^2&=&-xU_{yz}+zU_{xy}&=\R_5\R_8\\
\Q_{41}^2&=&U_{xy}&=\R_6\R_8\\
\Q_{42}^2&=&-xU_{zz}+zU_{xz}&=\R_5\R_9\\
\Q_{43}^2&=&-xU_{yz}+yU_{xz}&=\R_4\R_9\\
\Q_{44}^2&=&U_{xz}&=\R_6\R_9\\
\Q_{45}^2&=&y^2U_{zz}-2yzU_{yz}+z^2U_{yy}+2tU_{yy}
+2tU_{zz}&=\R_7\R_7+\R_3\R_9+\R_2\R_8\\

\Q_{46}^2&=&-yU_{yz}+zU_{yy}&=\R_7\R_8\\
\Q_{47}^2&=&U_{yy}&=\R_8\R_8\\
\Q_{48}^2&=&-yU_{zz}+zU_{yz}&=\R_7\R_9\\
\Q_{49}^2&=&U_{yz}&=\R_8\R_9\\
\Q_{50}^2&=&U_{zz}&=\R_9\R_9\\
\end{array}
\en

\subsection{Third order}

We let the characteristic depending up to third order

$Q(x,y,z,t,U,U_x,U_y,U_z,U_{xx},U_{xy},U_{xz},U_{yy},U_{yz},
,U_{zz},U_{xxx},U_{xxy},U_{xxz},U_{xyy},U_{xyz},
,U_{xzz},U_{yyy},U_{yyz},
,U_{yzz},U_{zzz})$

solving the invariance condition \eqref{inv}, we got the following set of characteristic
\bn
\begin{array}{llll}
\Q^3_1&=& 8t^3U_{xxx}+12t^2xU_{xx}+6tx^2U_x+x^3U+12t^2U_x+6txU\\

\Q^3_2&=&  4t^3U_{xxy}+4t^2xU_{xy}+2t^2yU_{xx}+tx^2U_y+2txyU_x+\( \frac{x^2y}{2} +ty\)U+2t^2U_y    \\

\Q^3_3&=& 4t^3U_{xxz}+4t^2xU_{xz}+2t^2zU_{xx}+tx^2U_z+2txzU_x 
+\( \frac{x^2z}{2}+tz \)U+2t^2U_z     \\
\Q^3_4&=&   4t^2xU_{xzz}-4t^2zU_{xxz}+2tx^2U_{zz}-2tz^2U_{xx}+x^2zU_z
-xz^2U_x-2t^2U_{xx}+4t^2U_{zz}+2tzU_z+\frac{x^2}{2}U+tU    \\

\Q^3_5&=&   -24t^3U_{xxy}+8t^3U_{yyy}-24t^2xU_{xy}-12t^2yU_{xx}+12t^2yU_{yy}
-6tx^2U_{y}-12txyU_{x}+  6ty^2U_y-3x^2yU+y^3U \\

\Q^3_6&=&  24t^3U_{xyy}-8t^3U_{xxx}+24t^2yU_{xy}-12t^2xU_{xx}+12t^2xU_{yy}
-6tx^2U_{x}+12txyU_{y}+  6ty^2U_x+3xy^2U-x^3U    \\

\Q^3_7&=& 4t^3(U_{yyz}-U_{xxz})+4t^2(yU_{yz}-xU_{xz})+2t^2z(U_{yy}-U_{xx
} )   +t(y^2-x^2)U_z-2txzU_x+2tyzU_y+\( \frac{y^2z}{2} -\frac{x^2 z}{2}\)U \\

\Q^3_8&=& 4t^2(yU_{yzz}-xU_{xzz})+4t^2z(U_{xxz}-U_{yyz})+2t(y^2-x^2)U_{zz}+2tz^2(U_{xx}-U_{yy})+xz(zU_x-xU_z)\\&&+yz(yU_z-zU_y)+2t^2(U_{xx}-U_{yy})+\( \frac{y^2 }{2}-\frac{x^2}{2}\)U \\

\Q^3_9&=& 4t^3(U_{yzz}-U_{xxy})-4t^2(xU_{xy}+yU_{xx})+2t^2yU_{zz}+4t^2zU_{yz}-tx^2U_y-2txyU_x+2tyzU_z+tz^2U_y+\\&&+\(\frac{yz^2}{2}-\frac{x^2y}{2}\)U\\

\Q^3_{10}&=&Uxyz+8t^3U_{xyz}+4t^2(xU_{yz}+yU_{xz}+zU_{xy})+
2txyU_z+2txzU_y+2tyzU_x \\

\Q^3_{11}&=& 4t^2yU_{xxy}-4t^2zU_{xxy}+4txyU_{xz}-4txzU_{xy}+x^2yU_z-x^2zU_y+Uyz+4t^2U_{yz}+4tyU_z\\

\Q^3_{12}&=&  4t^2xU_{xxy}-4t^2yU_{xxx}+16t^2yU_{xzz}-16t^2zU_{xyz}+4tx^2U_{xy}-4txyU
_{xx}+8txyU_{zz}-8txzU_{yz}+8tyzU_{xz}\\&&-8tz^2U_{xy}+x^3U_y-x^2yU_x+4xyzU_z-4xz^2U_y+Uxy-4t^2U_{xy}\\
\Q^3_{13}&=&4tx(yU_{xzz}-zU_{xyz})+4tz(zU_{xxy}-yU_{xxz})+2x^2(yU_{zz}-zU_{yz})+2xz(zU_{xy}-yU_{xz})+4t^2U_{xxy}+4tyU_{zz}\\&&-4tzU_{yz}-x^2U_y+Uy \\

\Q^3_{14}&=&-24t^3U_{xxz}+8t^3U_{zzz}-24t^2xU_{xz}-12t^2zU_{xx}+12t^2z
U_{zz}-6tx^2U_z-12txzU_x+6tz^2U_z+\( z^3-3x^2z \)U\\

\Q^3_{15}&=&-8t^3U_{xxx}+24t^3U_{xzz}-12t^2xU_{xx}+12t^2xU_{zz}+24t^2z
U_{xz}-6tx^2U_x+12txzU_z+6tz^2U_x+(3xz^2-x^3)U\\

\Q^3_{16}&=&-4t^2xU_{xzz}+4t^2zU_{xxz}-2tx^2U_{zz}+2tz^2U_{xx}-x^2zU_{z}+xz^2U_x+2t^2(U_{xx}-U_{zz})+\( \frac{z^2}{2}-\frac{x^2}{2}  \)U\\

\Q^3_{17}&=& 4t^2(xU_{xxz}-zU_{xxx})+4tx(xU_{xz}-zU_{xx})+x^3U_z-x^2zU_x+Uxz+12t^2U_{xz}+8txU_z\\

\Q^3_{18}&=& -4t^2U_{xxz}-4txU_{xz}-x^2U_z+Uz\\

\Q^3_{19}&=&4tx(xU_{xyy}-2yU_{xxy}+2yU_{yzz}-2zU_{yyz})+4ty^2(U_{xxx}-3
U_{xzz})+16tyzU_{xyz}-4tz^2U_{xyy}+2x^3U_{yy}-4x^2yU_{xy}\\&&+2xy^2(U_{xx}

-U_{zz})+8xyzU_{yz}-6xz^2U_{yy}-4y^2zU_{xz}+4yz^2U_{xy}+4t^2(U_{xxx}-
6U_{xzz})-4t(xU_{zz}+2zU_{xz})\\&&-x^2U_x+Ux\\

\Q^3_{20}&=&U\\

\Q^3_{21}&=&-2\,{t}^{2}xU_{{xxy}}+2\,{t}^{2}yU_{{xxx}}-6\,{t}^{2}yU_{{xzz}}+
6\,{t}^{2}zU_{{xyz}}-2\,t{x}^{2}U_{{xy}}+2\,txyU_{{xx}}-3\,txyU_{{
zz}}+3\,txzU_{{yz}}-3\,tyzU_{{xz}}\\&&+3\,t{z}^{2}U_{{xy}}-\frac{1}{2}\,{x}^{3
}U_{{y}}+\frac{1}{2}\,{x}^{2}yU_{{x}}-\frac{3}{2}\,xyzU_{{z}}+\frac{3}{2}\,x{z}^{2}U_{{y}}+2\,
{t}^{2}U_{{xy}}+tyU_{{x}}

 \\

\Q^3_{22}&=& -2\,{t}^{2}xU_{{x,x,z}}+2\,{t}^{2}zU_{{x,x,x}}-2\,t{x}^{2}U_{{x,z}}+2
\,txzU_{{x,x}}-\frac{1}{2}{x}^{3}U_{{z}}+\frac{1}{2}{x}^{2}zU_{{x}}-4\,{t}^{2}U_{{
x,z}}-3\,txU_{{z}}+tzU_{{x}}
\\

\Q_{23}^3 &=& -2\,{t}^{2}xU_{{xzz}}+2\,{t}^{2}zU_{{xxz}}-t{x}^{2}U_{{zz}}+t{z}^
{2}U_{{xx}}-\frac{1}{2}\,{x}^{2}zU_{{z}}+\frac{1}{2}\,x{z}^{2}U_{{x}}+2\,{t}^{2}U_{{x
x}}-2\,{t}^{2}U_{{zz}}+txU_{{x}}-tzU_{{z}}
\\
\end{array}
\en
\bn
\begin{array}{llllll}

\Q_{24}^3 &=&-2\,t{x}^{2}U_{{xyy}}+4\,txyU_{{xxy}}-4\,txyU_{{yzz}}+4\,txzU_{{
yyz}}-2\,t{y}^{2}U_{{xxx}}+6\,t{y}^{2}U_{{xzz}}-8\,tyzU_{{xyz}
}+2\,t{z}^{2}U_{{xyy}}-{x}^{3}U_{{yy}}\\&&+2\,{x}^{2}yU_{{xy}}-x{y}^{2
}U_{{xx}}+x{y}^{2}U_{{zz}}-4\,xyzU_{{yz}}+3\,x{z}^{2}U_{{yy}}+2\,{
y}^{2}zU_{{xz}}-2\,y{z}^{2}U_{{xy}}-2\,{t}^{2}U_{{xxx}}+12\,{t}^{2
}U_{{xzz}}\\&&+2\,txU_{{zz}}+4\,tzU_{{xz}}+\frac{1}{2}{x}^{2}U_{{x}}+tU_{{x}
}
\\

\Q_{25}^3 &=& 2\,{t}^{2}xU_{{xxy}}-\frac{2}{3}\,{t}^{2}xU_{{yyy}}-2\,{t}^{2}yU_{{xxx}}
+\frac{2}{3}\,{t}^{2}yU_{{xyy}}+\frac{16}{3}\,{t}^{2}yU_{{xzz}}-\frac{16}{3}{t}^{2}zU_{{
xyz}}+2\,t{x}^{2}U_{{xy}}-2\,txyU_{{xx}}-\frac{2}{3}\,txyU_{{yy}}\\&&+\frac{8}{3}\,tx
yU_{{zz}}-\frac{8}{3}\,txzU_{{yz}}+\frac{2}{3}\,t{y}^{2}U_{{xy}}+\frac{8}{3}\,tyzU_{{xz}}-
\frac{8}{3}t{z}^{2}U_{{xy}}+\frac{1}{2}\,{x}^{3}U_{{y}}-
\frac{1}{2}{x}^{2}yU_{{x}}-\frac{1}{6}x
{y}^{2}U_{{y}}+\frac{3}{4}xyzU_{{z}}-\frac{3}{4}x{z}^{2}U_{{y}}\\&&+\frac{1}{6}{y}^{3}U_{{x}
}
\\

\Q_{26}^3 &=& -2\,{t}^{2}xU_{{xyy}}+2\,{t}^{2}xU_{{xzz}}+2\,{t}^{2}yU_{{xxy}}-
2\,{t}^{2}yU_{{yzz}}-2\,{t}^{2}zU_{{xxz}}+2\,{t}^{2}zU_{{yyz}}-t
{x}^{2}U_{{yy}}+t{x}^{2}U_{{zz}}+t{y}^{2}U_{{xx}}-\\&&t{y}^{2}U_{{zz}}
-t{z}^{2}U_{{xx}}+t{z}^{2}U_{{yy}}-\frac{1}{2}\,{x}^{2}yU_{{y}}+\frac{1}{2}{x}^{2}
zU_{{z}}+\frac{1}{2}x{y}^{2}U_{{x}}-\frac{1}{2}\,x{z}^{2}U_{{x}}-\frac{1}{2}\,{y}^{2}zU_{{z}
}+\frac{1}{2}\,y{z}^{2}U_{{y}}
\\

\Q_{27}^3&=&  2\,{t}^{2}xU_{{xxz}}-2\,{t}^{2}xU_{{yyz}}-2\,{t}^{2}zU_{{xxx}}+2
\,{t}^{2}zU_{{xyy}}+2\,t{x}^{2}U_{{xz}}-2\,txyU_{{yz}}-2\,txzU_{{x
,x}}+2\,tyzU_{{xy}}+\frac{1}{2}\,{x}^{3}U_{{z}}-\\&&\frac{1}{2}\,{x}^{2}zU_{{x}}-\frac{1}{2}\,x{y
}^{2}U_{{z}}+\frac{1}{2}\,{y}^{2}zU_{{x}}+4\,{t}^{2}U_{{xz}}+2\,txU_{{z}}
\\

\Q_{28}^3 &=& 2\,t{x}^{2}U_{{xyy}}-4\,txyU_{{xxy}}+2\,txyU_{{yzz}}-2\,txzU_{{y
yz}}+2\,t{y}^{2}U_{{xxx}}-6\,t{y}^{2}U_{{xzz}}+10\,tyzU_{{xyz}
}-4\,t{z}^{2}U_{{xyy}}+{x}^{3}U_{{yy}}\\&&-2\,{x}^{2}yU_{{xy}}+x{y}^{2
}U_{{xx}}-2\,x{y}^{2}U_{{zz}}+5\,xyzU_{{yz}}-3\,x{z}^{2}U_{{yy}}-{
y}^{2}zU_{{xz}}+y{z}^{2}U_{{xy}}+2\,{t}^{2}U_{{xxx}}-2\,{t}^{2}U_{
{xyy}}-\\&&12\,{t}^{2}U_{{xzz}}-4\,txU_{{zz}}-2\,tzU_{{xz}}-\frac{1}{2}\,{x}
^{2}U_{{x}}+\frac{1}{2}\,{y}^{2}U_{{x}}
\\

\Q_{29}^3&=&2\,{t}^{2}xU_{{xxy}}-2\,{t}^{2}xU_{{yzz}}-2\,{t}^{2}yU_{{xxx}}+6
\,{t}^{2}yU_{{xzz}}-4\,{t}^{2}zU_{{xyz}}+2\,t{x}^{2}U_{{xy}}-2\,t
xyU_{{xx}}+2\,txyU_{{zz}}-4\,txzU_{{yz}}\\&&+4\,tyzU_{{xz}}-2\,t{z}^{2
}U_{{xy}}+\frac{1}{2}\,{x}^{3}U_{{y}}-\frac{1}{2}\,{x}^{2}yU_{{x}}+xyzU_{{z}}-\frac{3}{2}\,x{
z}^{2}U_{{y}}+\frac{3}{2}y{z}^{2}U_{{x}}
\\

\Q_{30}^3 &=& -4\,{t}^{2}xU_{{xyz}}+4\,{t}^{2}zU_{{xxy}}-2\,t{x}^{2}U_{{yz}}-2
\,txyU_{{xz}}+2\,txzU_{{xy}}+2\,tyzU_{{xx}}-{x}^{2}yU_{{z}}+xyzU_{{
x}}-4\,{t}^{2}U_{{yz}}-2\,tyU_{{z}} \\

\Q_{31}^3&= & -2\,txyU_{{yyz}}+2\,txzU_{{yyy}}+2\,t{y}^{2}U_{{xyz}}-2\,tyzU_{{
xyy}}-x{y}^{2}U_{{yz}}+xyzU_{{yy}}+{y}^{3}U_{{xz}}-{y}^{2}zU_{{x
y}}+8\,{t}^{2}U_{{xyz}}-\\&&2\,txU_{{yz}}+8\,tyU_{{xz}}+yzU_{{x}}
\\
\Q_{32}^3&=& -2\,t{x}^{2}U_{{{\it yyy}}}+4\,txyU_{{{\it xyy}}}-10\,txyU_{{{\it xzz}
}}+10\,txzU_{{{\it xyz}}}-2\,t{y}^{2}U_{{{\it xxy}}}+2\,t{y}^{2}U_{{{
\it yzz}}}+10\,tyzU_{{{\it xxz}}}-4\,tyzU_{{{\it yyz}}}-\\&&10\,t{z}^{2}U_
{{{\it xxy}}}+2\,t{z}^{2}U_{{{\it yyy}}}-{x}^{2}yU_{{{\it yy}}}-5\,{x}
^{2}yU_{{{\it zz}}}+5\,{x}^{2}zU_{{{\it yz}}}+2\,x{y}^{2}U_{{{\it xy}}
}+5\,xyzU_{{{\it xz}}}-5\,x{z}^{2}U_{{{\it xy}}}-{y}^{3}U_{{{\it xx}}}
\\&&+{y}^{3}U_{{{\it zz}}}-2\,{y}^{2}zU_{{{\it yz}}}+y{z}^{2}U_{{{\it yy}}
}-12\,{t}^{2}U_{{{\it xxy}}}+4\,{t}^{2}U_{{{\it yzz}}}-4\,tyU_{{{\it 
zz}}}+6\,tzU_{{{\it yz}}}+3\,xyU_{{x}}
\\

\Q_{33}^3&=& -x{y}^{2}U_{{{\it yzz}}}+2\,xyzU_{{{\it yyz}}}-x{z}^{2}U_{{{\it yyy}}}
+{y}^{3}U_{{{\it xzz}}}-2\,{y}^{2}zU_{{{\it xyz}}}+y{z}^{2}U_{{{\it 
xyy}}}-2\,txU_{{{\it yzz}}}+6\,tyU_{{{\it xzz}}}-4\,tzU_{{{\it xyz}}}\\&&+
xyU_{{{\it yy}}}-{y}^{2}U_{{{\it xy}}}+yU_{{x}}

\\

\Q_{34}^3&=& 12\,{t}^{2}xU_{{{\it xxz}}}-4\,{t}^{2}xU_{{{\it zzz}}}-12\,{t}^{2}zU_{
{{\it xxx}}}+4\,{t}^{2}zU_{{{\it xzz}}}+12\,t{x}^{2}U_{{{\it xz}}}-12
\,txzU_{{{\it xx}}}-4\,txzU_{{{\it zz}}}+4\,t{z}^{2}U_{{{\it xz}}}+3\,
{x}^{3}U_{{z}}\\&&-3\,{x}^{2}zU_{{x}}-x{z}^{2}U_{{z}}+{z}^{3}U_{{x}}+32\,{
t}^{2}U_{{{\it xz}}}+16\,txU_{{z}}
\\

\Q_{35}^3&=&4\,t{x}^{2}U_{{{\it xyy}}}-8\,txyU_{{{\it xxy}}}+8\,txyU_{{{\it yzz}}}
-8\,txzU_{{{\it yyz}}}+4\,t{y}^{2}U_{{{\it xxx}}}-12\,t{y}^{2}U_{{{
\it xzz}}}+16\,tyzU_{{{\it xyz}}}-4\,t{z}^{2}U_{{{\it xyy}}}+2\,{x}^{3
}U_{{{\it yy}}}\\&&-4\,{x}^{2}yU_{{{\it xy}}}+2\,x{y}^{2}U_{{{\it xx}}}-2
\,x{y}^{2}U_{{{\it zz}}}+8\,xyzU_{{{\it yz}}}-6\,x{z}^{2}U_{{{\it yy}}
}-4\,{y}^{2}zU_{{{\it xz}}}+4\,y{z}^{2}U_{{{\it xy}}}+4\,{t}^{2}U_{{{
\it xxx}}}-20\,{t}^{2}U_{{{\it xzz}}}\\&&-4\,txU_{{{\it zz}}}-4\,tzU_{{{
\it xz}}}-{x}^{2}U_{{x}}+{z}^{2}U_{{x}}

\\

\Q_{36}^3&=&-6\,t{x}^{2}U_{{{\it yyz}}}+12\,txzU_{{{\it xyy}}}+6\,t{y}^{2}U_{{{
\it xxz}}}-2\,t{y}^{2}U_{{{\it zzz}}}-12\,tyzU_{{{\it xxy}}}+4\,tyzU_{
{{\it yzz}}}-2\,t{z}^{2}U_{{{\it yyz}}}-6\,{x}^{2}yU_{{{\it yz}}}+3\,{
x}^{2}zU_{{{\it yy}}}\\&&+6\,x{y}^{2}U_{{{\it xz}}}-3\,{y}^{2}zU_{{{\it xx
}}}-{y}^{2}zU_{{{\it zz}}}+2\,y{z}^{2}U_{{{\it yz}}}-{z}^{3}U_{{{\it 
yy}}}+24\,{t}^{2}U_{{{\it xxz}}}-16\,{t}^{2}U_{{{\it yyz}}}-4\,{t}^{2}
U_{{{\it zzz}}}+18\,txU_{{{\it xz}}}-8\,tyU_{{{\it yz}}}\\&&-2\,tzU_{{{
\it zz}}}+3\,xzU_{{x}}
\\
\Q_{37}^3&=&  2\,txU_{{{\it yyz}}}-2\,tyU_{{{\it xyz}}}+xyU_{{{\it yz}}}-{y}^{2}U_{{
{\it xz}}}+zU_{{x}}
   \\

\end{array}
\en

\bn
\begin{array}{lllll}
\Q_{38}^3&=& -2\,{x}^{2}yU_{{{\it yzz}}}+2\,{x}^{2}zU_{{{\it yyz}}}+2\,x{y}^{2}U_{{
{\it xzz}}}-2\,x{z}^{2}U_{{{\it xyy}}}-2\,{y}^{2}zU_{{{\it xxz}}}+2\,y
{z}^{2}U_{{{\it xxy}}}+4\,txU_{{{\it xzz}}}-4\,tzU_{{{\it xxz}}}\\&&+{x}^{
2}U_{{{\it yy}}}-{y}^{2}U_{{{\it xx}}}+{y}^{2}U_{{{\it zz}}}-{z}^{2}U_
{{{\it yy}}}+2\,tU_{{{\it zz}}}+xU_{{x}}
    \\

\Q_{39}^3&=&   U_x  \\

\Q_{40}^3&=&   -4\,{t}^{2}yU_{{{\it yzz}}}+4\,{t}^{2}zU_{{{\it yyz}}}-2\,t{y}^{2}U_{{
{\it zz}}}+2\,t{z}^{2}U_{{{\it yy}}}-{y}^{2}zU_{{z}}+y{z}^{2}U_{{y}}+4
\,{t}^{2}U_{{{\it yy}}}-4\,{t}^{2}U_{{{\it zz}}}+2\,tyU_{{y}}-2\,tzU_{
{z}}
  \\

\Q_{41}^3&=&   -4\,{t}^{2}yU_{{{\it xxz}}}+4\,{t}^{2}zU_{{{\it xxy}}}-4\,txyU_{{{\it 
xz}}}+4\,txzU_{{{\it xy}}}-{x}^{2}yU_{{z}}+{x}^{2}zU_{{y}}-2\,tyU_{{z}
}+2\,tzU_{{y}}
  \\

\Q_{42}^3&=&   -4\,{t}^{2}yU_{{{\it xzz}}}+4\,{t}^{2}zU_{{{\it xyz}}}-2\,txyU_{{{\it 
zz}}}+2\,txzU_{{{\it yz}}}-2\,tyzU_{{{\it xz}}}+2\,t{z}^{2}U_{{{\it xy
}}}-xyzU_{{z}}+x{z}^{2}U_{{y}}+4\,{t}^{2}U_{{{\it xy}}}+2\,txU_{{y}}
  \\

\Q_{43}^3&=&  -4\,txyU_{{{\it xzz}}}+4\,txzU_{{{\it xyz}}}+4\,tyzU_{{{\it xxz}}}-4\,
t{z}^{2}U_{{{\it xxy}}}-2\,{x}^{2}yU_{{{\it zz}}}+2\,{x}^{2}zU_{{{\it 
yz}}}+2\,xyzU_{{{\it xz}}}-2\,x{z}^{2}U_{{{\it xy}}}\\&&-4\,{t}^{2}U_{{{
\it xxy}}}-4\,tyU_{{{\it zz}}}+4\,tzU_{{{\it yz}}}+{x}^{2}U_{{y}}+2\,t
U_{{y}}
   \\

\Q_{44}^3&=&   4\,{t}^{2}yU_{{{\it xxz}}}-4\,{t}^{2}yU_{{{\it yyz}}}-4\,{t}^{2}zU_{{{
\it xxy}}}+4\,{t}^{2}zU_{{{\it yyy}}}+4\,txyU_{{{\it xz}}}-4\,txzU_{{{
\it xy}}}-4\,t{y}^{2}U_{{{\it yz}}}+4\,tyzU_{{{\it yy}}}+{x}^{2}yU_{{z
}}\\&&-{x}^{2}zU_{{y}}-{y}^{3}U_{{z}}+{y}^{2}zU_{{y}}-8\,{t}^{2}U_{{{\it 
yz}}}-4\,tyU_{{z}}

  \\

\Q_{45}^3&=&  4\,txyU_{{{\it xzz}}}-4\,txzU_{{{\it xyz}}}-4\,t{y}^{2}U_{{{\it yzz}}}
-4\,tyzU_{{{\it xxz}}}+8\,tyzU_{{{\it yyz}}}+4\,t{z}^{2}U_{{{\it xxy}}
}-4\,t{z}^{2}U_{{{\it yyy}}}+2\,{x}^{2}yU_{{{\it zz}}}\\&&-2\,{x}^{2}zU_{{
{\it yz}}}-2\,xyzU_{{{\it xz}}}+2\,x{z}^{2}U_{{{\it xy}}}-2\,{y}^{3}U_
{{{\it zz}}}+4\,{y}^{2}zU_{{{\it yz}}}-2\,y{z}^{2}U_{{{\it yy}}}+4\,{t
}^{2}U_{{{\it xxy}}}-4\,{t}^{2}U_{{{\it yyy}}}-8\,{t}^{2}U_{{{\it yzz}
}}\\&&-8\,tyU_{{{\it zz}}}+4\,tzU_{{{\it yz}}}-{x}^{2}U_{{y}}+{y}^{2}U_{{y
}}
   \\

\Q_{46}^3&=&   -4\,{t}^{2}yU_{{{\it xyz}}}+4\,{t}^{2}zU_{{{\it xyy}}}-2\,txyU_{{{\it 
yz}}}+2\,txzU_{{{\it yy}}}-2\,t{y}^{2}U_{{{\it xz}}}+2\,tyzU_{{{\it xy
}}}-x{y}^{2}U_{{z}}+xyzU_{{y}}-4\,{t}^{2}U_{{{\it xz}}}-2\,txU_{{z}}
  \\

\Q_{47}^3&=& -2\,txyU_{{{\it yzz}}}+2\,txzU_{{{\it yyz}}}-2\,t{y}^{2}U_{{{\it xzz}}
}+6\,tyzU_{{{\it xyz}}}-4\,t{z}^{2}U_{{{\it xyy}}}-2\,x{y}^{2}U_{{{
\it zz}}}+3\,xyzU_{{{\it yz}}}-x{z}^{2}U_{{{\it yy}}}+{y}^{2}zU_{{{
\it xz}}}\\&&-y{z}^{2}U_{{{\it xy}}}-4\,{t}^{2}U_{{{\it xyy}}}-4\,{t}^{2}U
_{{{\it xzz}}}-4\,txU_{{{\it zz}}}+2\,tzU_{{{\it xz}}}+xyU_{{y}}
   \\

\Q_{48}^3&=&  -2\,t{y}^{2}U_{{{\it zzz}}}+4\,tyzU_{{{\it yzz}}}-2\,t{z}^{2}U_{{{\it 
yyz}}}-{y}^{2}zU_{{{\it zz}}}+2\,y{z}^{2}U_{{{\it yz}}}-{z}^{3}U_{{{
\it yy}}}+8\,{t}^{2}U_{{{\it yyz}}}-4\,{t}^{2}U_{{{\it zzz}}}+10\,tyU_
{{{\it yz}}}\\&&-2\,tzU_{{{\it zz}}}+3\,yzU_{{y}}
   \\

\Q_{49}^3&=& 4\,tyU_{{{\it yzz}}}-4\,tzU_{{{\it yyz}}}+{y}^{2}U_{{{\it zz}}}-{z}^{2
}U_{{{\it yy}}}+2\,tU_{{{\it zz}}}+yU_{{y}}
    \\

\Q_{50}^3&=&  12\,{t}^{2}yU_{{{\it xxz}}}-4\,{t}^{2}yU_{{{\it zzz}}}-12\,{t}^{2}zU_{
{{\it xxy}}}+4\,{t}^{2}zU_{{{\it yzz}}}+12\,txyU_{{{\it xz}}}-12\,txzU
_{{{\it xy}}}-4\,tyzU_{{{\it zz}}}+4\,t{z}^{2}U_{{{\it yz}}}+3\,{x}^{2
}yU_{{z}}\\&&-3\,{x}^{2}zU_{{y}}-y{z}^{2}U_{{z}}+{z}^{3}U_{{y}}+8\,{t}^{2}
U_{{{\it yz}}}+4\,tyU_{{z}}
   \\

\Q_{51}^3&=&  4\,txyU_{{{\it xzz}}}-4\,txzU_{{{\it xyz}}}-4\,tyzU_{{{\it xxz}}}+4\,t
{z}^{2}U_{{{\it xxy}}}+2\,{x}^{2}yU_{{{\it zz}}}-2\,{x}^{2}zU_{{{\it 
yz}}}-2\,xyzU_{{{\it xz}}}+2\,x{z}^{2}U_{{{\it xy}}}+4\,{t}^{2}U_{{{
\it xxy}}}\\&&+4\,{t}^{2}U_{{{\it yzz}}}+4\,tyU_{{{\it zz}}}-{x}^{2}U_{{y}
}+{z}^{2}U_{{y}}
   \\

\Q_{52}^3&=&  -2\,txyU_{{{\it yyz}}}+2\,txzU_{{{\it yyy}}}+2\,t{y}^{2}U_{{{\it xyz}}
}-2\,tyzU_{{{\it xyy}}}-x{y}^{2}U_{{{\it yz}}}+xyzU_{{{\it yy}}}+{y}^{
3}U_{{{\it xz}}}-{y}^{2}zU_{{{\it xy}}}+8\,{t}^{2}U_{{{\it xyz}}}+\\&&6\,t
yU_{{{\it xz}}}+xzU_{{y}}

   \\

\Q_{53}^3&=&  2\,tyU_{{{\it zzz}}}-2\,tzU_{{{\it yzz}}}+yzU_{{{\it zz}}}-{z}^{2}U_{{
{\it yz}}}+zU_{{y}}
   \\

\Q_{54}^3&=&   -x{y}^{2}U_{{{\it yzz}}}+2\,xyzU_{{{\it yyz}}}-x{z}^{2}U_{{{\it yyy}}}
+{y}^{3}U_{{{\it xzz}}}-2\,{y}^{2}zU_{{{\it xyz}}}+y{z}^{2}U_{{{\it 
xyy}}}-2\,txU_{{{\it yzz}}}+6\,tyU_{{{\it xzz}}}-4\,tzU_{{{\it xyz}}}\\&&+
xyU_{{{\it yy}}}-{y}^{2}U_{{{\it xy}}}+xU_{{y}}
  \\
  
  \Q_{55}^3&=&   U_y  \\

\Q_{56}^3&=&  4\,{t}^{2}U_{{{\it xxz}}}+4\,txU_{{{\it xz}}}+{x}^{2}U_{{z}}+2\,tU_{{z
}}

   \\

\Q_{57}^3&=&  -4\,{t}^{2}U_{{{\it xxz}}}+4\,{t}^{2}U_{{{\it zzz}}}-4\,txU_{{{\it xz}
}}+4\,tzU_{{{\it zz}}}-{x}^{2}U_{{z}}+{z}^{2}U_{{z}}

   \\
\end{array}
\en

\bn
\begin{array}{llll}
\Q_{58}^3&=&   4\,{t}^{2}U_{{{\it yzz}}}+2\,tyU_{{{\it zz}}}+2\,tzU_{{{\it yz}}}+yzU_
{{z}}

  \\

\Q_{59}^3&=&  4\,{t}^{2}U_{{{\it xzz}}}+2\,txU_{{{\it zz}}}+2\,tzU_{{{\it xz}}}+xzU_
{{z}}
   \\

\Q_{60}^3&=&   -2\,{t}^{2}U_{{{\it xxz}}}+2\,{t}^{2}U_{{{\it yyz}}}-2\,txU_{{{\it xz}
}}+2\,tyU_{{{\it yz}}}-1/2\,{x}^{2}U_{{z}}+1/2\,{y}^{2}U_{{z}}+4\,tU_{
{{\it zz}}}+2\,zU_{{z}}
  \\

\Q_{61}^3&=& -4\,{t}^{2}U_{{{\it xxz}}}+4\,{t}^{2}U_{{{\it yyz}}}-4\,txU_{{{\it xz}
}}+4\,tyU_{{{\it yz}}}-{x}^{2}U_{{z}}+{y}^{2}U_{{z}}

   \\

\Q_{62}^3&=&  4\,{t}^{2}U_{{{\it xyz}}}+2\,txU_{{{\it yz}}}+2\,tyU_{{{\it xz}}}+xyU_
{{z}}
   \\

\Q_{63}^3&=& 2\,tyU_{{{\it zzz}}}-2\,tzU_{{{\it yzz}}}+yzU_{{{\it zz}}}-{z}^{2}U_{{
{\it yz}}}+yU_{{z}}
    \\

\Q_{64}^3&=&   2\,txU_{{{\it yyz}}}-2\,tyU_{{{\it xyz}}}+xyU_{{{\it yz}}}-{y}^{2}U_{{
{\it xz}}}+xU_{{z}}
  \\

\Q_{65}^3&=&  U_z   \\

\Q_{66}^3&=&  2\,t{x}^{2}U_{{{\it xyy}}}-4\,txyU_{{{\it xxy}}}+4\,txyU_{{{\it yzz}}}
-4\,txzU_{{{\it yyz}}}+2\,t{y}^{2}U_{{{\it xxx}}}-6\,t{y}^{2}U_{{{\it 
xzz}}}+8\,tyzU_{{{\it xyz}}}-2\,t{z}^{2}U_{{{\it xyy}}}\\&&+{x}^{3}U_{{{
\it yy}}}-2\,{x}^{2}yU_{{{\it xy}}}+x{y}^{2}U_{{{\it xx}}}-x{y}^{2}U_{
{{\it zz}}}+4\,xyzU_{{{\it yz}}}-3\,x{z}^{2}U_{{{\it yy}}}-2\,{y}^{2}z
U_{{{\it xz}}}+2\,y{z}^{2}U_{{{\it xy}}}+4\,{t}^{2}U_{{{\it xxx}}}\\&&-12
\,{t}^{2}U_{{{\it xzz}}}+2\,txU_{{{\it xx}}}-2\,txU_{{{\it zz}}}-4\,tz
U_{{{\it xz}}}
   \\

\Q_{67}^3&=&  -2\,t{x}^{2}U_{{{\it xyy}}}+2\,t{x}^{2}U_{{{\it xzz}}}+4\,txyU_{{{\it 
xxy}}}-4\,txyU_{{{\it yzz}}}-4\,txzU_{{{\it xxz}}}+4\,txzU_{{{\it yyz}
}}-2\,t{y}^{2}U_{{{\it xxx}}}+6\,t{y}^{2}U_{{{\it xzz}}}\\&&-8\,tyzU_{{{
\it xyz}}}+2\,t{z}^{2}U_{{{\it xxx}}}+2\,t{z}^{2}U_{{{\it xyy}}}-{x}^{
3}U_{{{\it yy}}}+{x}^{3}U_{{{\it zz}}}+2\,{x}^{2}yU_{{{\it xy}}}-2\,{x
}^{2}zU_{{{\it xz}}}-x{y}^{2}U_{{{\it xx}}}+x{y}^{2}U_{{{\it zz}}}\\&&-4\,
xyzU_{{{\it yz}}}+x{z}^{2}U_{{{\it xx}}}+3\,x{z}^{2}U_{{{\it yy}}}+2\,
{y}^{2}zU_{{{\it xz}}}-2\,y{z}^{2}U_{{{\it xy}}}+16\,{t}^{2}U_{{{\it 
xzz}}}+8\,txU_{{{\it zz}}}

   \\

\Q_{68}^3&=&  2\,t{x}^{2}U_{{{\it xyz}}}-2\,txyU_{{{\it xxz}}}+2\,txyU_{{{\it yyz}}}
-2\,txzU_{{{\it xxy}}}-2\,txzU_{{{\it yyy}}}-2\,t{y}^{2}U_{{{\it xyz}}
}+2\,tyzU_{{{\it xxx}}}+2\,tyzU_{{{\it xyy}}}\\&&+{x}^{3}U_{{{\it yz}}}-{x
}^{2}yU_{{{\it xz}}}-{x}^{2}zU_{{{\it xy}}}+x{y}^{2}U_{{{\it yz}}}+xyz
U_{{{\it xx}}}-xyzU_{{{\it yy}}}-{y}^{3}U_{{{\it xz}}}+{y}^{2}zU_{{{
\it xy}}}+8\,txU_{{{\it yz}}}-8\,tyU_{{{\it xz}}}
   \\

\Q_{69}^3&=&   3\,{x}^{3}U_{{{\it yyz}}}-6\,{x}^{2}yU_{{{\it xyz}}}-3\,{x}^{2}zU_{{{
\it xyy}}}+3\,x{y}^{2}U_{{{\it xxz}}}-x{y}^{2}U_{{{\it zzz}}}+6\,xyzU_
{{{\it xxy}}}+2\,xyzU_{{{\it yzz}}}-x{z}^{2}U_{{{\it yyz}}}\\&&-3\,{y}^{2}
zU_{{{\it xxx}}}+{y}^{2}zU_{{{\it xzz}}}-2\,y{z}^{2}U_{{{\it xyz}}}+{z
}^{3}U_{{{\it xyy}}}+10\,txU_{{{\it yyz}}}-2\,txU_{{{\it zzz}}}-10\,ty
U_{{{\it xyz}}}+2\,tzU_{{{\it xzz}}}\\&&-3\,{x}^{2}U_{{{\it xz}}}-3\,xyU_{
{{\it yz}}}+3\,xzU_{{{\it xx}}}+3\,{y}^{2}U_{{{\it xz}}}
  \\

\Q_{70}^3&=&  {x}^{3}U_{{{\it yyy}}}-3\,{x}^{2}yU_{{{\it xyy}}}+3\,x{y}^{2}U_{{{\it 
xxy}}}-{y}^{3}U_{{{\it xxx}}}-3\,{x}^{2}U_{{{\it xy}}}+3\,xyU_{{{\it 
xx}}}-3\,xyU_{{{\it yy}}}+3\,{y}^{2}U_{{{\it xy}}}

   \\ 

\Q_{71}^3&=&  -{x}^{2}U_{{{\it xyy}}}+2\,xyU_{{{\it xxy}}}+2\,xyU_{{{\it yzz}}}-2\,x
zU_{{{\it yyz}}}-{y}^{2}U_{{{\it xxx}}}-{y}^{2}U_{{{\it xzz}}}+{z}^{2}
U_{{{\it xyy}}}-2\,tU_{{{\it xzz}}}+xU_{{{\it xx}}}

   \\ 

\Q_{72}^3&=& 2\,t{x}^{2}U_{{{\it yyy}}}-4\,txyU_{{{\it xyy}}}+4\,txyU_{{{\it xzz}}}
-4\,txzU_{{{\it xyz}}}+2\,t{y}^{2}U_{{{\it xxy}}}-2\,t{y}^{2}U_{{{\it 
yzz}}}-4\,tyzU_{{{\it xxz}}}+4\,tyzU_{{{\it yyz}}}\\&&+4\,t{z}^{2}U_{{{
\it xxy}}}-2\,t{z}^{2}U_{{{\it yyy}}}+{x}^{2}yU_{{{\it yy}}}+2\,{x}^{2
}yU_{{{\it zz}}}-2\,{x}^{2}zU_{{{\it yz}}}-2\,x{y}^{2}U_{{{\it xy}}}-2
\,xyzU_{{{\it xz}}}+2\,x{z}^{2}U_{{{\it xy}}}+{y}^{3}U_{{{\it xx}}}\\&&-{y
}^{3}U_{{{\it zz}}}+2\,{y}^{2}zU_{{{\it yz}}}-y{z}^{2}U_{{{\it yy}}}+
12\,{t}^{2}U_{{{\it xxy}}}-4\,{t}^{2}U_{{{\it yzz}}}+6\,tyU_{{{\it xx}
}}-2\,tyU_{{{\it zz}}}
    \\ 

\Q_{73}^3&=&  6\,t{x}^{2}U_{{{\it yyz}}}-12\,txzU_{{{\it xyy}}}-6\,t{y}^{2}U_{{{\it 
xxz}}}+2\,t{y}^{2}U_{{{\it zzz}}}+12\,tyzU_{{{\it xxy}}}-4\,tyzU_{{{
\it yzz}}}+2\,t{z}^{2}U_{{{\it yyz}}}+6\,{x}^{2}yU_{{{\it yz}}}\\&&-3\,{x}
^{2}zU_{{{\it yy}}}-6\,x{y}^{2}U_{{{\it xz}}}+3\,{y}^{2}zU_{{{\it xx}}
}+{y}^{2}zU_{{{\it zz}}}-2\,y{z}^{2}U_{{{\it yz}}}+{z}^{3}U_{{{\it yy}
}}-12\,{t}^{2}U_{{{\it xxz}}}+16\,{t}^{2}U_{{{\it yyz}}}+4\,{t}^{2}U_{
{{\it zzz}}}\\&&-12\,txU_{{{\it xz}}}+8\,tyU_{{{\it yz}}}+6\,tzU_{{{\it xx
}}}+2\,tzU_{{{\it zz}}}

   \\ 
   
\Q_{74}^3&=&  2\,{x}^{2}yU_{{{\it yzz}}}-2\,{x}^{2}zU_{{{\it yyz}}}-2\,x{y}^{2}U_{{{
\it xzz}}}+2\,x{z}^{2}U_{{{\it xyy}}}+2\,{y}^{2}zU_{{{\it xxz}}}-2\,y{
z}^{2}U_{{{\it xxy}}}-4\,txU_{{{\it xzz}}}+4\,tzU_{{{\it xxz}}}\\&&-{x}^{2
}U_{{{\it yy}}}+{y}^{2}U_{{{\it xx}}}-{y}^{2}U_{{{\it zz}}}+{z}^{2}U_{
{{\it yy}}}+2\,tU_{{{\it xx}}}-2\,tU_{{{\it zz}}}
   \\

\Q_{75}^3&=&   -6\,t{x}^{2}U_{{{\it yyz}}}+2\,t{x}^{2}U_{{{\it zzz}}}+12\,txzU_{{{
\it xyy}}}-4\,txzU_{{{\it xzz}}}+6\,t{y}^{2}U_{{{\it xxz}}}-2\,t{y}^{2
}U_{{{\it zzz}}}-12\,tyzU_{{{\it xxy}}}+4\,tyzU_{{{\it yzz}}}\\&&+2\,t{z}^
{2}U_{{{\it xxz}}}-2\,t{z}^{2}U_{{{\it yyz}}}-6\,{x}^{2}yU_{{{\it yz}}
}+3\,{x}^{2}zU_{{{\it yy}}}+{x}^{2}zU_{{{\it zz}}}+6\,x{y}^{2}U_{{{
\it xz}}}-2\,x{z}^{2}U_{{{\it xz}}}-3\,{y}^{2}zU_{{{\it xx}}}-{y}^{2}z
U_{{{\it zz}}}\\&&+2\,y{z}^{2}U_{{{\it yz}}}+{z}^{3}U_{{{\it xx}}}-{z}^{3}
U_{{{\it yy}}}+16\,{t}^{2}U_{{{\it xxz}}}-16\,{t}^{2}U_{{{\it yyz}}}+8
\,txU_{{{\it xz}}}-8\,tyU_{{{\it yz}}}

  \\ 
\end{array}
\en

\bn
\begin{array}{llll}
\Q_{76}^3&=&  -2\,t{x}^{2}U_{{{\it yyy}}}+6\,t{x}^{2}U_{{{\it yzz}}}+4\,txyU_{{{\it 
xyy}}}-4\,txyU_{{{\it xzz}}}-8\,txzU_{{{\it xyz}}}-2\,t{y}^{2}U_{{{
\it xxy}}}+2\,t{y}^{2}U_{{{\it yzz}}}+4\,tyzU_{{{\it xxz}}}\\&&-4\,tyzU_{{
{\it yyz}}}+2\,t{z}^{2}U_{{{\it xxy}}}+2\,t{z}^{2}U_{{{\it yyy}}}-{x}^
{2}yU_{{{\it yy}}}+{x}^{2}yU_{{{\it zz}}}+2\,{x}^{2}zU_{{{\it yz}}}+2
\,x{y}^{2}U_{{{\it xy}}}-4\,xyzU_{{{\it xz}}}\\&&-2\,x{z}^{2}U_{{{\it xy}}
}-{y}^{3}U_{{{\it xx}}}+{y}^{3}U_{{{\it zz}}}-2\,{y}^{2}zU_{{{\it yz}}
}+3\,y{z}^{2}U_{{{\it xx}}}+y{z}^{2}U_{{{\it yy}}}+16\,{t}^{2}U_{{{
\it yzz}}}+8\,tyU_{{{\it zz}}}
   \\ 

\Q_{77}^3&=&  -2\,{x}^{2}yU_{{{\it yzz}}}+2\,{x}^{2}zU_{{{\it yyz}}}+2\,x{y}^{2}U_{{
{\it xzz}}}-2\,x{z}^{2}U_{{{\it xyy}}}-2\,{y}^{2}zU_{{{\it xxz}}}+2\,y
{z}^{2}U_{{{\it xxy}}}+{x}^{2}U_{{{\it yy}}}-{x}^{2}U_{{{\it zz}}}-{y}
^{2}U_{{{\it xx}}}\\&&+{y}^{2}U_{{{\it zz}}}+{z}^{2}U_{{{\it xx}}}-{z}^{2}
U_{{{\it yy}}}

   \\ 

\Q_{78}^3&=&   3\,{x}^{2}yU_{{{\it yyz}}}-3\,{x}^{2}zU_{{{\it yyy}}}-6\,x{y}^{2}U_{{{
\it xyz}}}+6\,xyzU_{{{\it xyy}}}+3\,{y}^{3}U_{{{\it xxz}}}-{y}^{3}U_{{
{\it zzz}}}-3\,{y}^{2}zU_{{{\it xxy}}}+3\,{y}^{2}zU_{{{\it yzz}}}\\&&-3\,y
{z}^{2}U_{{{\it yyz}}}+{z}^{3}U_{{{\it yyy}}}-24\,txU_{{{\it xyz}}}+24
\,tyU_{{{\it xxz}}}-12\,tyU_{{{\it zzz}}}+12\,tzU_{{{\it yzz}}}-3\,{x}
^{2}U_{{{\it yz}}}+3\,yzU_{{{\it xx}}}-\\&&3\,yzU_{{{\it zz}}}+3\,{z}^{2}U
_{{{\it yz}}}

  \\ 

\Q_{79}^3&=&  -{x}^{2}U_{{{\it yyz}}}+2\,xyU_{{{\it xyz}}}-{y}^{2}U_{{{\it xxz}}}+{y
}^{2}U_{{{\it zzz}}}-2\,yzU_{{{\it yzz}}}+{z}^{2}U_{{{\it yyz}}}+2\,tU
_{{{\it zzz}}}+zU_{{{\it xx}}}
   \\ 

\Q_{80}^3&=&   -{x}^{2}U_{{{\it yyy}}}+2\,xyU_{{{\it xyy}}}-{y}^{2}U_{{{\it xxy}}}+{y
}^{2}U_{{{\it yzz}}}-2\,yzU_{{{\it yyz}}}+{z}^{2}U_{{{\it yyy}}}+2\,tU
_{{{\it yzz}}}+yU_{{{\it xx}}}

  \\      
 
\Q_{81}^3&=&U_{xx}\\
\Q_{82}^3&=&  2\,txyU_{{{\it yzz}}}-2\,txzU_{{{\it yyz}}}-2\,tyzU_{{{\it xyz}}}+2\,t
{z}^{2}U_{{{\it xyy}}}+x{y}^{2}U_{{{\it zz}}}-xyzU_{{{\it yz}}}-{y}^{2
}zU_{{{\it xz}}}+y{z}^{2}U_{{{\it xy}}}+4\,{t}^{2}U_{{{\it xyy}}}\\&&+2\,t
xU_{{{\it zz}}}+2\,tyU_{{{\it xy}}}-2\,tzU_{{{\it xz}}}
 \\

\Q_{83}^3&=& 2\,txyU_{{{\it xzz}}}-2\,txzU_{{{\it xyz}}}-2\,tyzU_{{{\it xxz}}}+2\,t
{z}^{2}U_{{{\it xxy}}}+{x}^{2}yU_{{{\it zz}}}-{x}^{2}zU_{{{\it yz}}}-x
yzU_{{{\it xz}}}+x{z}^{2}U_{{{\it xy}}}+4\,{t}^{2}U_{{{\it xxy}}}\\&&+2\,t
xU_{{{\it xy}}}+2\,tyU_{{{\it zz}}}-2\,tzU_{{{\it yz}}}
   \\

\Q_{84}^3&=&  2\,txyU_{{{\it yyz}}}-2\,txzU_{{{\it yyy}}}-2\,t{y}^{2}U_{{{\it xyz}}}
+2\,tyzU_{{{\it xyy}}}+x{y}^{2}U_{{{\it yz}}}-xyzU_{{{\it yy}}}-{y}^{3
}U_{{{\it xz}}}+{y}^{2}zU_{{{\it xy}}}-4\,{t}^{2}U_{{{\it xyz}}}\\&&+2\,tx
U_{{{\it yz}}}-6\,tyU_{{{\it xz}}}+2\,tzU_{{{\it xy}}}
  \\

\Q_{85}^3&=&  x{y}^{2}U_{{{\it yzz}}}-2\,xyzU_{{{\it yyz}}}+x{z}^{2}U_{{{\it yyy}}}-
{y}^{3}U_{{{\it xzz}}}+2\,{y}^{2}zU_{{{\it xyz}}}-y{z}^{2}U_{{{\it xyy
}}}+2\,txU_{{{\it yzz}}}-6\,tyU_{{{\it xzz}}}+4\,tzU_{{{\it xyz}}}\\&&-xyU
_{{{\it yy}}}+{y}^{2}U_{{{\it xy}}}+2\,tU_{{{\it xy}}}
  \\

\Q_{86}^3&=&  6\,txyU_{{{\it xyz}}}-6\,txzU_{{{\it xyy}}}-6\,t{y}^{2}U_{{{\it xxz}}}
+2\,t{y}^{2}U_{{{\it zzz}}}+6\,tyzU_{{{\it xxy}}}-4\,tyzU_{{{\it yzz}}
}+2\,t{z}^{2}U_{{{\it yyz}}}+3\,{x}^{2}yU_{{{\it yz}}}\\&&-3\,{x}^{2}zU_{{
{\it yy}}}-3\,x{y}^{2}U_{{{\it xz}}}+3\,xyzU_{{{\it xy}}}+{y}^{2}zU_{{
{\it zz}}}-2\,y{z}^{2}U_{{{\it yz}}}+{z}^{3}U_{{{\it yy}}}-12\,{t}^{2}
U_{{{\it xxz}}}+4\,{t}^{2}U_{{{\it yyz}}}+4\,{t}^{2}U_{{{\it zzz}}}\\&&-6
\,txU_{{{\it xz}}}+2\,tyU_{{{\it yz}}}+2\,tzU_{{{\it zz}}}
  \\

\Q_{87}^3&=&  {x}^{2}yU_{{{\it yzz}}}-{x}^{2}zU_{{{\it yyz}}}-x{y}^{2}U_{{{\it xzz}}
}+x{z}^{2}U_{{{\it xyy}}}+{y}^{2}zU_{{{\it xxz}}}-y{z}^{2}U_{{{\it xxy
}}}-2\,txU_{{{\it xzz}}}-2\,tyU_{{{\it yzz}}}+2\,tzU_{{{\it xxz}}}\\&&+2\,
tzU_{{{\it yyz}}}-{x}^{2}U_{{{\it yy}}}+xyU_{{{\it xy}}}-{y}^{2}U_{{{
\it zz}}}+{z}^{2}U_{{{\it yy}}}-2\,tU_{{{\it zz}}}
  \\

\Q_{88}^3&=&  3\,{x}^{2}yU_{{{\it yyz}}}-3\,{x}^{2}zU_{{{\it yyy}}}-6\,x{y}^{2}U_{{{
\it xyz}}}+6\,xyzU_{{{\it xyy}}}+3\,{y}^{3}U_{{{\it xxz}}}-{y}^{3}U_{{
{\it zzz}}}-3\,{y}^{2}zU_{{{\it xxy}}}+3\,{y}^{2}zU_{{{\it yzz}}}\\&&-3\,y
{z}^{2}U_{{{\it yyz}}}+{z}^{3}U_{{{\it yyy}}}-18\,txU_{{{\it xyz}}}+18
\,tyU_{{{\it xxz}}}-12\,tyU_{{{\it zzz}}}+12\,tzU_{{{\it yzz}}}-3\,{x}
^{2}U_{{{\it yz}}}+3\,xzU_{{{\it xy}}}\\&&-3\,yzU_{{{\it zz}}}+3\,{z}^{2}U
_{{{\it yz}}}
  \\

\Q_{89}^3&=&  -{x}^{2}U_{{{\it yyy}}}+2\,xyU_{{{\it xyy}}}-{y}^{2}U_{{{\it xxy}}}+{y
}^{2}U_{{{\it yzz}}}-2\,yzU_{{{\it yyz}}}+{z}^{2}U_{{{\it yyy}}}+2\,tU
_{{{\it yzz}}}+xU_{{{\it xy}}}
  \\
  
  \Q_{90}^3&=&  -6\,txyU_{{{\it yyz}}}+2\,txyU_{{{\it zzz}}}+6\,txzU_{{{\it yyy}}}-2\,
txzU_{{{\it yzz}}}+6\,t{y}^{2}U_{{{\it xyz}}}-6\,tyzU_{{{\it xyy}}}-2
\,tyzU_{{{\it xzz}}}+2\,t{z}^{2}U_{{{\it xyz}}}\\&&-3\,x{y}^{2}U_{{{\it yz
}}}+3\,xyzU_{{{\it yy}}}+xyzU_{{{\it zz}}}-x{z}^{2}U_{{{\it yz}}}+3\,{
y}^{3}U_{{{\it xz}}}-3\,{y}^{2}zU_{{{\it xy}}}-y{z}^{2}U_{{{\it xz}}}+
{z}^{3}U_{{{\it xy}}}+16\,{t}^{2}U_{{{\it xyz}}}\\&&-8\,txU_{{{\it yz}}}+
16\,tyU_{{{\it xz}}}
  \\

\Q_{91}^3&=&   -x{y}^{2}U_{{{\it yzz}}}+2\,xyzU_{{{\it yyz}}}-x{z}^{2}U_{{{\it yyy}}}
+{y}^{3}U_{{{\it xzz}}}-2\,{y}^{2}zU_{{{\it xyz}}}+y{z}^{2}U_{{{\it 
xyy}}}-4\,txU_{{{\it yzz}}}+4\,tyU_{{{\it xzz}}}+xyU_{{{\it yy}}}\\&&-xyU_
{{{\it zz}}}-{y}^{2}U_{{{\it xy}}}+{z}^{2}U_{{{\it xy}}}
 \\
\end{array}
\en

\bn
\begin{array}{llllll}

\Q_{92}^3&=&  -2\,xyzU_{{{\it yzz}}}+x{z}^{2}U_{{{\it yyz}}}-{y}^{2}zU_{{{\it xzz}}}
+2\,y{z}^{2}U_{{{\it xyz}}}-{z}^{3}U_{{{\it xyy}}}-4\,txU_{{{\it yyz}}
}+2\,txU_{{{\it zzz}}}+4\,tyU_{{{\it xyz}}}-2\,tzU_{{{\it xzz}}}\\&&-3\,xy
U_{{{\it yz}}}+{y}^{2}U_{{{\it zzz}}}+3\,yzU_{{{\it xy}}}
  \\

\Q_{93}^3&=&  -xyU_{{{\it zzz}}}+xzU_{{{\it yzz}}}+yzU_{{{\it xzz}}}-{z}^{2}U_{{{
\it xyz}}}+zU_{{{\it xy}}}
  \\

\Q_{94}^3&=&  -2\,xyU_{{{\it yzz}}}+2\,xzU_{{{\it yyz}}}+{y}^{2}U_{{{\it xzz}}}-{z}^
{2}U_{{{\it xyy}}}+2\,tU_{{{\it xzz}}}+yU_{{{\it xy}}}
  \\

\Q_{95}^3&=&  U_{xy}  \\

\Q_{96}^3&=&  -2\,txU_{{{\it xzz}}}+2\,tzU_{{{\it xxz}}}-{x}^{2}U_{{{\it zz}}}+xzU_{
{{\it xz}}}-2\,tU_{{{\it zz}}}
  \\

\Q_{97}^3&=& -2\,txU_{{{\it xyz}}}+2\,tyU_{{{\it xxz}}}-2\,tyU_{{{\it zzz}}}+2\,tzU
_{{{\it yzz}}}-{x}^{2}U_{{{\it yz}}}+xyU_{{{\it xz}}}-yzU_{{{\it zz}}}
+{z}^{2}U_{{{\it yz}}}
   \\

\Q_{98}^3&=&   -{x}^{2}U_{{{\it yyz}}}+2\,xyU_{{{\it xyz}}}-{y}^{2}U_{{{\it xxz}}}+{y
}^{2}U_{{{\it zzz}}}-2\,yzU_{{{\it yzz}}}+{z}^{2}U_{{{\it yyz}}}+2\,tU
_{{{\it zzz}}}+xU_{{{\it xz}}}
 \\

\Q_{99}^3&=&   -2\,txU_{{{\it yyz}}}+2\,tyU_{{{\it xyz}}}-xyU_{{{\it yz}}}+{y}^{2}U_{
{{\it xz}}}+2\,tU_{{{\it xz}}}
 \\

\Q_{100}^3&=&  2\,txU_{{{\it yyz}}}-2\,txU_{{{\it zzz}}}-2\,tyU_{{{\it xyz}}}+2\,tzU_
{{{\it xzz}}}+xyU_{{{\it yz}}}-xzU_{{{\it zz}}}-{y}^{2}U_{{{\it xz}}}+
{z}^{2}U_{{{\it xz}}}
  \\

\Q_{101}^3&=& -2\,txU_{{{\it yzz}}}+2\,tzU_{{{\it xyz}}}-xyU_{{{\it zz}}}+yzU_{{{
\it xz}}}
  \\

\Q_{102}^3&=&   2\,tU_{{{\it xzz}}}+zU_{{{\it xz}}} \\

\Q_{103}^3&=&  -xyU_{{{\it zzz}}}+xzU_{{{\it yzz}}}+yzU_{{{\it xzz}}}-{z}^{2}U_{{{
\it xyz}}}+yU_{{{\it xz}}}

  \\

\Q_{104}^3&=&  U_{xz}  \\

\Q_{105}^3&=&  2\,t{y}^{2}U_{{{\it yzz}}}-4\,tyzU_{{{\it yyz}}}+2\,t{z}^{2}U_{{{\it 
yyy}}}+{y}^{3}U_{{{\it zz}}}-2\,{y}^{2}zU_{{{\it yz}}}+y{z}^{2}U_{{{
\it yy}}}+4\,{t}^{2}U_{{{\it yyy}}}+4\,{t}^{2}U_{{{\it yzz}}}+2\,tyU_{
{{\it yy}}}\\&&+6\,tyU_{{{\it zz}}}-4\,tzU_{{{\it yz}}}

  \\

\Q_{106}^3&=&  {y}^{3}U_{{{\it zzz}}}-3\,{y}^{2}zU_{{{\it yzz}}}+3\,y{z}^{2}U_{{{\it 
yyz}}}-{z}^{3}U_{{{\it yyy}}}-3\,{y}^{2}U_{{{\it yz}}}+3\,yzU_{{{\it 
yy}}}-3\,yzU_{{{\it zz}}}+3\,{z}^{2}U_{{{\it yz}}}
  \\

\Q_{107}^3&=&  -{y}^{2}U_{{{\it yzz}}}+2\,yzU_{{{\it yyz}}}-{z}^{2}U_{{{\it yyy}}}-2
\,tU_{{{\it yzz}}}+yU_{{{\it yy}}}
  \\

\Q_{108}^3&=&   2\,t{y}^{2}U_{{{\it xzz}}}-4\,tyzU_{{{\it xyz}}}+2\,t{z}^{2}U_{{{\it 
xyy}}}+x{y}^{2}U_{{{\it zz}}}-2\,xyzU_{{{\it yz}}}+x{z}^{2}U_{{{\it yy
}}}+4\,{t}^{2}U_{{{\it xyy}}}+4\,{t}^{2}U_{{{\it xzz}}}+2\,txU_{{{\it 
yy}}}\\&&+2\,txU_{{{\it zz}}}
 \\

\Q_{109}^3&=&  2\,t{y}^{2}U_{{{\it zzz}}}-4\,tyzU_{{{\it yzz}}}+2\,t{z}^{2}U_{{{\it 
yyz}}}+{y}^{2}zU_{{{\it zz}}}-2\,y{z}^{2}U_{{{\it yz}}}+{z}^{3}U_{{{
\it yy}}}+4\,{t}^{2}U_{{{\it yyz}}}+4\,{t}^{2}U_{{{\it zzz}}}-4\,tyU_{
{{\it yz}}}\\&&+6\,tzU_{{{\it yy}}}+2\,tzU_{{{\it zz}}}
  \\

\Q_{110}^3&=& -4\,tyU_{{{\it yzz}}}+4\,tzU_{{{\it yyz}}}-{y}^{2}U_{{{\it zz}}}+{z}^{
2}U_{{{\it yy}}}+2\,tU_{{{\it yy}}}-2\,tU_{{{\it zz}}}
\\
\Q_{111}^3&=&  x{y}^{2}U_{{{\it zzz}}}-2\,xyzU_{{{\it yzz}}}+x{z}^{2}U_{{{\it yyz}}}-
{y}^{2}zU_{{{\it xzz}}}+2\,y{z}^{2}U_{{{\it xyz}}}-{z}^{3}U_{{{\it xyy
}}}+2\,txU_{{{\it yyz}}}+2\,txU_{{{\it zzz}}}-2\,tyU_{{{\it xyz}}}\\&&-2\,
tzU_{{{\it xzz}}}-3\,xyU_{{{\it yz}}}+3\,xzU_{{{\it yy}}}
  \\

\Q_{112}^3&=& -2\,xyU_{{{\it yzz}}}+2\,xzU_{{{\it yyz}}}+{y}^{2}U_{{{\it xzz}}}-{z}^
{2}U_{{{\it xyy}}}+2\,tU_{{{\it xzz}}}+xU_{{{\it yy}}}
   \\

\Q_{113}^3&=&   -{y}^{2}U_{{{\it zzz}}}+2\,yzU_{{{\it yzz}}}-{z}^{2}U_{{{\it yyz}}}-2
\,tU_{{{\it zzz}}}+zU_{{{\it yy}}}
 \\

\Q_{114}^3&=&   U_{yy} \\

\Q_{115}^3&=& -2\,tyU_{{{\it zzz}}}+2\,tzU_{{{\it yzz}}}-yzU_{{{\it zz}}}+{z}^{2}U_{
{{\it yz}}}+2\,tU_{{{\it yz}}}
   \\

\Q_{116}^3&=&  -2\,tyU_{{{\it yzz}}}+2\,tzU_{{{\it yyz}}}-{y}^{2}U_{{{\it zz}}}+yzU_{
{{\it yz}}}-2\,tU_{{{\it zz}}}
  \\

\Q_{117}^3&=&  -{y}^{2}U_{{{\it zzz}}}+2\,yzU_{{{\it yzz}}}-{z}^{2}U_{{{\it yyz}}}-2
\,tU_{{{\it zzz}}}+yU_{{{\it yz}}}
  \\
  
  \Q_{118}^3&=&  -2\,tyU_{{{\it xzz}}}+2\,tzU_{{{\it xyz}}}-xyU_{{{\it zz}}}+xzU_{{{
\it yz}}}
  \\

\Q_{119}^3&=&  2\,tU_{{{\it yzz}}}+zU_{{{\it yz}}}  \\

\end{array}
\en

\bn
\begin{array}{lllll}
\Q_{120}^3&=&  -xyU_{{{\it zzz}}}+xzU_{{{\it yzz}}}+yzU_{{{\it xzz}}}-{z}^{2}U_{{{
\it xyz}}}+xU_{{{\it yz}}}
   \\
\Q_{121}^3&=&  U_{yz}   \\

\Q_{122}^3&=&  2\,tU_{{{\it zzz}}}+zU_{{{\it zz}}}   \\

\Q_{123}^3&=&  2\,tU_{{{\it xzz}}}+xU_{{{\it zz}}}   \\

\Q_{124}^3&=&   2\,tU_{{{\it yzz}}}+yU_{{{\it zz}}}  \\

\Q_{125}^3&=&  U_{zz}   \\

\Q_{126}^3&=&   -3\,{x}^{3}U_{{{\it yyz}}}+6\,{x}^{2}yU_{{{\it xyz}}}+3\,{x}^{2}zU_{{{
\it xyy}}}-3\,x{y}^{2}U_{{{\it xxz}}}+x{y}^{2}U_{{{\it zzz}}}-6\,xyzU_
{{{\it xxy}}}-2\,xyzU_{{{\it yzz}}}+x{z}^{2}U_{{{\it yyz}}}\\&&+3\,{y}^{2}
zU_{{{\it xxx}}}-{y}^{2}zU_{{{\it xzz}}}+2\,y{z}^{2}U_{{{\it xyz}}}-{z
}^{3}U_{{{\it xyy}}}-6\,txU_{{{\it xxz}}}-16\,txU_{{{\it yyz}}}+2\,txU
_{{{\it zzz}}}+16\,tyU_{{{\it xyz}}}+\\&&6\,tzU_{{{\it xxx}}}-2\,tzU_{{{
\it xzz}}}
  \\

\Q_{127}^3&=&  -{x}^{3}U_{{{\it yyy}}}+3\,{x}^{2}yU_{{{\it xyy}}}-3\,x{y}^{2}U_{{{
\it xxy}}}+3\,x{y}^{2}U_{{{\it yzz}}}-6\,xyzU_{{{\it yyz}}}+3\,x{z}^{2
}U_{{{\it yyy}}}+{y}^{3}U_{{{\it xxx}}}-3\,{y}^{3}U_{{{\it xzz}}}\\&&+6\,{
y}^{2}zU_{{{\it xyz}}}-3\,y{z}^{2}U_{{{\it xyy}}}-6\,txU_{{{\it xxy}}}
+6\,txU_{{{\it yzz}}}+6\,tyU_{{{\it xxx}}}-18\,tyU_{{{\it xzz}}}+12\,t
zU_{{{\it xyz}}}
   \\

\Q_{128}^3&=&  {x}^{2}U_{{{\it xyy}}}-2\,xyU_{{{\it xxy}}}-2\,xyU_{{{\it yzz}}}+2\,xz
U_{{{\it yyz}}}+{y}^{2}U_{{{\it xxx}}}+{y}^{2}U_{{{\it xzz}}}-{z}^{2}U
_{{{\it xyy}}}+2\,tU_{{{\it xxx}}}+2\,tU_{{{\it xzz}}}
   \\

\Q_{129}^3&=&  3\,{x}^{3}U_{{{\it yyz}}}-{x}^{3}U_{{{\it zzz}}}-6\,{x}^{2}yU_{{{\it 
xyz}}}-3\,{x}^{2}zU_{{{\it xyy}}}+3\,{x}^{2}zU_{{{\it xzz}}}+3\,x{y}^{
2}U_{{{\it xxz}}}-x{y}^{2}U_{{{\it zzz}}}+6\,xyzU_{{{\it xxy}}}\\&&+2\,xyz
U_{{{\it yzz}}}-3\,x{z}^{2}U_{{{\it xxz}}}-x{z}^{2}U_{{{\it yyz}}}-3\,
{y}^{2}zU_{{{\it xxx}}}+{y}^{2}zU_{{{\it xzz}}}-2\,y{z}^{2}U_{{{\it 
xyz}}}+{z}^{3}U_{{{\it xxx}}}+{z}^{3}U_{{{\it xyy}}}+\\&&16\,txU_{{{\it 
yyz}}}-8\,txU_{{{\it zzz}}}-16\,tyU_{{{\it xyz}}}+8\,tzU_{{{\it xzz}}}

   \\

\Q_{130}^3&=&   {x}^{3}U_{{{\it yyy}}}-3\,{x}^{3}U_{{{\it yzz}}}-3\,{x}^{2}yU_{{{\it 
xyy}}}+3\,{x}^{2}yU_{{{\it xzz}}}+6\,{x}^{2}zU_{{{\it xyz}}}+3\,x{y}^{
2}U_{{{\it xxy}}}-3\,x{y}^{2}U_{{{\it yzz}}}-6\,xyzU_{{{\it xxz}}}\\&&+6\,
xyzU_{{{\it yyz}}}-3\,x{z}^{2}U_{{{\it xxy}}}-3\,x{z}^{2}U_{{{\it yyy}
}}-{y}^{3}U_{{{\it xxx}}}+3\,{y}^{3}U_{{{\it xzz}}}-6\,{y}^{2}zU_{{{
\it xyz}}}+3\,y{z}^{2}U_{{{\it xxx}}}+3\,y{z}^{2}U_{{{\it xyy}}}\\&&-24\,t
xU_{{{\it yzz}}}+24\,tyU_{{{\it xzz}}}
  \\

\Q_{131}^3&=&   -{x}^{2}U_{{{\it xyy}}}+{x}^{2}U_{{{\it xzz}}}+2\,xyU_{{{\it xxy}}}+2
\,xyU_{{{\it yzz}}}-2\,xzU_{{{\it xxz}}}-2\,xzU_{{{\it yyz}}}-{y}^{2}U
_{{{\it xxx}}}-{y}^{2}U_{{{\it xzz}}}+{z}^{2}U_{{{\it xxx}}}\\&&+{z}^{2}U_
{{{\it xyy}}}
  \\

\Q_{132}^3&=&   {x}^{2}U_{{{\it xyz}}}-xyU_{{{\it xxz}}}-xyU_{{{\it zzz}}}-xzU_{{{\it 
xxy}}}+xzU_{{{\it yzz}}}+yzU_{{{\it xxx}}}+yzU_{{{\it xzz}}}-{z}^{2}U_
{{{\it xyz}}}
  \\

\Q_{133}^3&=&   -xU_{{{\it xxz}}}+zU_{{{\it xxx}}}  \\

\Q_{134}^3&=& -xU_{{{\it xxy}}}+yU_{{{\it xxx}}}    \\

\Q_{135}^3&=&   U_{xxx}  \\

\Q_{136}^3&=&  -{x}^{2}yU_{{{\it yzz}}}+{x}^{2}zU_{{{\it yyz}}}+x{y}^{2}U_{{{\it xzz}
}}-x{z}^{2}U_{{{\it xyy}}}-{y}^{2}zU_{{{\it xxz}}}+y{z}^{2}U_{{{\it 
xxy}}}-2\,txU_{{{\it xyy}}}+2\,txU_{{{\it xzz}}}+2\,tyU_{{{\it xxy}}}\\&&-
2\,tyU_{{{\it yzz}}}-2\,tzU_{{{\it xxz}}}+2\,tzU_{{{\it yyz}}}
   \\

\Q_{137}^3&=&   -3\,{x}^{2}yU_{{{\it yyz}}}+3\,{x}^{2}zU_{{{\it yyy}}}+6\,x{y}^{2}U_{{
{\it xyz}}}-6\,xyzU_{{{\it xyy}}}-3\,{y}^{3}U_{{{\it xxz}}}+{y}^{3}U_{
{{\it zzz}}}+3\,{y}^{2}zU_{{{\it xxy}}}\\&&-3\,{y}^{2}zU_{{{\it yzz}}}+3\,
y{z}^{2}U_{{{\it yyz}}}-{z}^{3}U_{{{\it yyy}}}+12\,txU_{{{\it xyz}}}-
18\,tyU_{{{\it xxz}}}+6\,tyU_{{{\it zzz}}}+6\,tzU_{{{\it xxy}}}-6\,tzU
_{{{\it yzz}}}
  \\

\Q_{138}^3&=&   {x}^{2}U_{{{\it yyy}}}-2\,xyU_{{{\it xyy}}}+{y}^{2}U_{{{\it xxy}}}-{y}
^{2}U_{{{\it yzz}}}+2\,yzU_{{{\it yyz}}}-{z}^{2}U_{{{\it yyy}}}+2\,tU_
{{{\it xxy}}}-2\,tU_{{{\it yzz}}}

  \\

\Q_{139}^3&=&   3\,{x}^{2}yU_{{{\it yyz}}}-{x}^{2}yU_{{{\it zzz}}}-3\,{x}^{2}zU_{{{
\it yyy}}}+{x}^{2}zU_{{{\it yzz}}}-6\,x{y}^{2}U_{{{\it xyz}}}+6\,xyzU_
{{{\it xyy}}}+2\,xyzU_{{{\it xzz}}}-2\,x{z}^{2}U_{{{\it xyz}}}\\&&+3\,{y}^
{3}U_{{{\it xxz}}}-{y}^{3}U_{{{\it zzz}}}-3\,{y}^{2}zU_{{{\it xxy}}}+3
\,{y}^{2}zU_{{{\it yzz}}}-y{z}^{2}U_{{{\it xxz}}}-3\,y{z}^{2}U_{{{\it 
yyz}}}+{z}^{3}U_{{{\it xxy}}}+{z}^{3}U_{{{\it yyy}}}-16\,txU_{{{\it 
xyz}}}\\&&+16\,tyU_{{{\it xxz}}}-8\,tyU_{{{\it zzz}}}+8\,tzU_{{{\it yzz}}}
  \\

\Q_{140}^3&=&  -{x}^{2}U_{{{\it yyy}}}+{x}^{2}U_{{{\it yzz}}}+2\,xyU_{{{\it xyy}}}-2
\,xzU_{{{\it xyz}}}-{y}^{2}U_{{{\it xxy}}}+{y}^{2}U_{{{\it yzz}}}-2\,y
zU_{{{\it yyz}}}+{z}^{2}U_{{{\it xxy}}}+{z}^{2}U_{{{\it yyy}}}+4\,tU_{
{{\it yzz}}}
 \\

\end{array}
\en

\bn
\begin{array}{lllll}

\Q_{141}^3&=& {x}^{2}U_{{{\it yyz}}}-xyU_{{{\it xyz}}}-xzU_{{{\it xyy}}}-{y}^{2}U_{{
{\it zzz}}}+yzU_{{{\it xxy}}}+2\,yzU_{{{\it yzz}}}-{z}^{2}U_{{{\it yyz
}}}-2\,tU_{{{\it zzz}}}
    \\

\Q_{142}^3&=&   -xU_{{{\it xyz}}}+zU_{{{\it xxy}}}  \\

\Q_{143}^3&=&  -xU_{{{\it xyy}}}+yU_{{{\it xxy}}}   \\

\Q_{144}^3&=&   U_{xxy}  \\

\Q_{145}^3&=&  {x}^{2}U_{{{\it yyz}}}-2\,xyU_{{{\it xyz}}}+{y}^{2}U_{{{\it xxz}}}-{y}
^{2}U_{{{\it zzz}}}+2\,yzU_{{{\it yzz}}}-{z}^{2}U_{{{\it yyz}}}+2\,tU_
{{{\it xxz}}}-2\,tU_{{{\it zzz}}}
  \\

\Q_{146}^3&=&   -{x}^{2}U_{{{\it yyz}}}+{x}^{2}U_{{{\it zzz}}}+2\,xyU_{{{\it xyz}}}-2
\,xzU_{{{\it xzz}}}-{y}^{2}U_{{{\it xxz}}}+{y}^{2}U_{{{\it zzz}}}-2\,y
zU_{{{\it yzz}}}+{z}^{2}U_{{{\it xxz}}}+{z}^{2}U_{{{\it yyz}}}+4\,tU_{
{{\it zzz}}}
  \\

\Q_{147}^3&=&  {x}^{2}U_{{{\it yzz}}}-xyU_{{{\it xzz}}}-xzU_{{{\it xyz}}}+yzU_{{{\it 
xxz}}}+2\,tU_{{{\it yzz}}}
   \\

\Q_{148}^3&=&   -xU_{{{\it xzz}}}+zU_{{{\it xxz}}}  \\

\Q_{149}^3&=&  -xU_{{{\it xyz}}}+yU_{{{\it xxz}}}   \\

\Q_{150}^3&=&    U_{xxz} \\

\Q_{151}^3&=&   -x{y}^{2}U_{{{\it yzz}}}+2\,xyzU_{{{\it yyz}}}-x{z}^{2}U_{{{\it yyy}}}
+{y}^{3}U_{{{\it xzz}}}-2\,{y}^{2}zU_{{{\it xyz}}}+y{z}^{2}U_{{{\it 
xyy}}}-2\,txU_{{{\it yyy}}}-2\,txU_{{{\it yzz}}}+2\,tyU_{{{\it xyy}}}\\&&+
6\,tyU_{{{\it xzz}}}-4\,tzU_{{{\it xyz}}}
  \\

\Q_{152}^3&=& xyU_{{{\it yyz}}}+xyU_{{{\it zzz}}}-xzU_{{{\it yyy}}}-xzU_{{{\it yzz}}
}-{y}^{2}U_{{{\it xyz}}}+yzU_{{{\it xyy}}}-yzU_{{{\it xzz}}}+{z}^{2}U_
{{{\it xyz}}}
    \\

\Q_{153}^3&=&  -xU_{{{\it yyy}}}+yU_{{{\it xyy}}}   \\

\Q_{154}^3&=&  -x{y}^{2}U_{{{\it zzz}}}+2\,xyzU_{{{\it yzz}}}-x{z}^{2}U_{{{\it yyz}}}
+{y}^{2}zU_{{{\it xzz}}}-2\,y{z}^{2}U_{{{\it xyz}}}+{z}^{3}U_{{{\it 
xyy}}}-2\,txU_{{{\it yyz}}}-2\,txU_{{{\it zzz}}}-4\,tyU_{{{\it xyz}}}\\&&+
6\,tzU_{{{\it xyy}}}+2\,tzU_{{{\it xzz}}}
   \\

\Q_{155}^3&=&   2\,xyU_{{{\it yzz}}}-2\,xzU_{{{\it yyz}}}-{y}^{2}U_{{{\it xzz}}}+{z}^{
2}U_{{{\it xyy}}}+2\,tU_{{{\it xyy}}}-2\,tU_{{{\it xzz}}}

  \\

\Q_{156}^3&=&  -xU_{{{\it yyz}}}+zU_{{{\it xyy}}}   \\

\Q_{157}^3&=&  U_{xyy}   \\

\Q_{158}^3&=&  xyU_{{{\it zzz}}}-xzU_{{{\it yzz}}}-yzU_{{{\it xzz}}}+{z}^{2}U_{{{\it 
xyz}}}+2\,tU_{{{\it xyz}}}
   \\

\Q_{159}^3&=&   xyU_{{{\it yzz}}}-xzU_{{{\it yyz}}}-{y}^{2}U_{{{\it xzz}}}+yzU_{{{\it 
xyz}}}-2\,tU_{{{\it xzz}}}
  \\

\Q_{160}^3&=&  -xU_{{{\it yyz}}}+yU_{{{\it xyz}}}   \\

\Q_{161}^3&=&    -xU_{{{\it yzz}}}+U_{{{\it xyz}}}z \\

\Q_{162}^3&=&   U_{xyz}  \\

\Q_{163}^3&=& -xU_{{{\it zzz}}}+zU_{{{\it xzz}}}    \\

\Q_{164}^3&=&   -xU_{{{\it yzz}}}+yU_{{{\it xzz}}}  \\

\Q_{165}^3&=&    U_{xzz} \\

\Q_{166}^3&=&   -{y}^{3}U_{{{\it zzz}}}+3\,{y}^{2}zU_{{{\it yzz}}}-3\,y{z}^{2}U_{{{
\it yyz}}}+{z}^{3}U_{{{\it yyy}}}-6\,tyU_{{{\it yyz}}}-6\,tyU_{{{\it 
zzz}}}+6\,tzU_{{{\it yyy}}}+6\,tzU_{{{\it yzz}}}
  \\

\Q_{167}^3&=&  {y}^{2}U_{{{\it yzz}}}-2\,yzU_{{{\it yyz}}}+{z}^{2}U_{{{\it yyy}}}+2\,
tU_{{{\it yyy}}}+2\,tU_{{{\it yzz}}}
   \\

\Q_{168}^3&=&  -yU_{{{\it yyz}}}+zU_{{{\it yyy}}}   \\

\Q_{169}^3&=&    U_{yyy} \\

\end{array}
\en

\bn
\begin{array}{llll}
\Q_{170}^3&=&   {y}^{2}U_{{{\it zzz}}}-2\,yzU_{{{\it yzz}}}+{z}^{2}U_{{{\it yyz}}}+2\,
tU_{{{\it yyz}}}+2\,tU_{{{\it zzz}}}
  \\

\Q_{171}^3&=&   -yU_{{{\it yzz}}}+zU_{{{\it yyz}}}  \\

\Q_{172}^3&=& U_{yyz}    \\

\Q_{173}^3&=&  -yU_{{{\it zzz}}}+zU_{{{\it yzz}}}   \\

\Q_{174}^3&=&   U_{yzz}  \\

\Q_{175}^3&=&  U_{zzz}   \\
\end{array}
\en

\section{Basis for conservation law}

The linear diffusion equation \eqref{heat}, viz. $U_t=U_{xx}+U_{yy}+U_{zz}$ has an obvious conserved vector $T$ with components namely
\bn T^1=-u, \quad T^2=u_x,\quad T^3=u_y,\quad T^4=u_z \label{consv}\en

We can generate an infinite set of conservation law associated to \eqref{heat}. We firstly determine the point symmetry associated with \eqref{consv}. To that end, we utilise the condition \eqref{tbr}. This yield, for $T^1,T^2,T^3$ and $T^4$ respectively, the following symmetries

\bn
\begin{array}{llll}
\X_4& = -x\p_y+y\p_x \\
\X_5& = -x\p_z +z\p_y \\
\X_6& = \p_x \\
\X_7& =  -y\p_z+z\p_y\\
\X_8& =  \p_y\\
\X_9& =  \p_z\\
\X_{11}& =\p_t  \\
\end{array}
\en
In fact \{$\X_4,X_5,\X_6,\X_7,\X_8,\X_9,\X_{11}$\} form a subalgebra of the Lie algebra \eqref{alge}. The remaining set of symmetries is \{$\X_1,\X_2,\X_3,\X_{10},\X_{12},\X_{13}$\} which can  be used to generate new non trivial conserved vectors. 

In \eqref{alge} we see that 

$[\X_1,\X_7]=[\X_2,\X_5]=[\X_3,\X_4]=[\X_4,\X_{10}]=[\X_4,\X_{12}]=[\X_4,\X_{13}]=0$

By virtue of theorem 1.1 none of them can generate a non trivial conserved vector.

The conserved vector $T$ with components given in \eqref{consv} is the only basis for the conservation law with respect to the group $G$ of transformations corresponding to \eqref{gen0}. In this equivalence class, we can generate infinite many conservation laws by applying successively the symmetries found in Section 2.

Eg. $U_{xx}\p_U$ is a second order generalized symmetry for the equation \eqref{heat}.

By \eqref{tbar},

$\bar{T}^1=-U_{xx},\quad \bar{T}^2=U_{xxx},\quad \bar{T}^3=U_{xxy},\quad \bar{T}^4=U_{xxz}$ are components of a conserved vector of \eqref{systm1}. In fact

$$D_t(-U_{xx})+D_x(U_{xxx})+D_y(U_{xxy})+D(U_{xxz})|_{U_t=U_{xx}+U_{yy}+U_{zz}}=0$$
\section{Conclusion}

Even if though, we were able to use the invariance condition to construct the basis of generalized symmetries in characteristic form  up to third order. The results matches with our predictions using the formula found in Section 2. We also showed that the linear diffusion equation \eqref{heat} has infinite many conservation laws in one equivalence class  with respect to the group $G$ of the equation \eqref{heat}.

\end{document}